\documentclass[11pt]{article}
\usepackage{exscale,makeidx,amssymb,amsmath,epsfig,epsf,psfrag,cases,multicol%bbm,color}
}
\usepackage{amsmath} % If you want to use AMS-LaTeX, comment out these
\usepackage{amssymb} % two commands.

\usepackage{amsthm}

\usepackage{graphics}
\usepackage{graphicx}
 \DeclareGraphicsExtensions{.eps, .ps}

\newtheorem{prop}{Proposition}
\newtheorem{defi}[prop]{Definition}
\newtheorem{lemm}[prop]{Lemma}
\newtheorem{theo}[prop]{Theorem}
\newtheorem{coro}[prop]{Corollary}
\newtheorem{rem}[prop]{Remark}

\newenvironment{pf}{\begin{trivlist}\item[] \textbf{Proof.}}
                     {\hspace*{\fill} $\square$\end{trivlist}}
\newenvironment{pfofprop4a}{\begin{trivlist}\item[] \textbf{Proof of Proposition \ref{prop4a}.}}
                     {\hspace*{\fill} $\square$\end{trivlist}}

\newcommand{\paramgamma}{h}

\newcommand{\beq}{\begin{eqnarray*}}
\newcommand{\eeq}{\end{eqnarray*}}

\newcommand{\bP}{\mathbb{P}}

\newcommand{\bN}{\mathbb{N}}
\newcommand{\bR}{\mathbb{R}}
\newcommand{\bD}{\mathbb{D}}

\newcommand{\bH}{\mathbb{H}}
\newcommand{\bE}{\mathbb{E}}

\newcommand{\bQ}{\mathbb{Q}}
\newcommand{\bT}{\mathbb{T}}

\newcommand{\bU}{\mathbb{U}}
\newcommand{\bM}{\mathbb{M}}

\newcommand{\cB}{\mathcal{B}}

\newcommand{\cF}{\mathcal{F}}
\newcommand{\cG}{\mathcal{G}}
\newcommand{\cH}{\mathcal{H}}

\newcommand{\cP}{\mathcal{P}}

\newcommand{\cR}{\mathcal{R}}
\newcommand{\cS}{\mathcal{S}}
\newcommand{\cT}{\mathcal{T}}

\newcommand{\cW}{\mathcal{W}}

\newcommand{\ft}{\mathbf{t}}

\begin{document}

\title{Growth of Galton-Watson trees: immigration and lifetimes}
\author{
Xiao'ou Cao\thanks{University of Oxford; email: cao@stats.ox.ac.uk}
\and %
Matthias Winkel\thanks{%
University of Oxford; email: winkel@stats.ox.ac.uk}}
\maketitle
\begin{abstract}

We study certain consistent families $(F_\lambda)_{\lambda\ge 0}$ of Galton-Watson forests with lifetimes as edge
lengths and/or immigrants as progenitors of the trees in $F_\lambda$. Specifically, consistency here refers to the property that for each $\mu\le\lambda$, the forest $F_\mu$ has the same distribution as the subforest of $F_\lambda$ spanned by the black leaves in a Bernoulli leaf colouring, where each leaf of $F_\lambda$ is coloured in black independently with probability $\mu/\lambda$. The case of exponentially distributed lifetimes and no immigration was studied by Duquesne and Winkel and related to the genealogy of Markovian continuous-state branching processes. We characterise here such families in the framework of arbitrary lifetime distributions and immigration according to a renewal process, related to Sagitov's (non-Markovian) generalisation of continuous-state branching renewal processes, and similar processes with immigration.

\emph{AMS 2000 subject classifications: 60J80.\newline
}

%{\tt xxx Draft between authors of work in progress. Not for general distribution. }\vspace{-0.7cm}
\end{abstract}

%\tableofcontents
\section{Introduction}

Galton-Watson branching processes are a classical model for the evolution of population sizes, see e.g. \cite{AtN,Har-55}.
More specifically, there is an interest in the underlying genealogical trees. In the most basic
model, there is a single \em progenitor \em that produces $i$ children with probability $q(i)$ for
some \em offspring distribution \em $q$ on $\bN_0=\{0,1,2,\ldots\}$; recursively, each individual
in the population produces children independently and according to the same distribution $q$. We
represent this model by a \em graph-theoretic tree \em rooted at the progenitor, where each individual is a vertex and the \em parent-child relation \em specifies edges $v\rightarrow w$ between parent $v$ and child $w$. Vertices related to just their parent vertex and to no child
vertices are called
\em leaves\em. More precisely, we will follow
Neveu \cite{Nev-GW} to distinguish individuals (see Section \ref{trees}). We
will consider in this paper the following well-known and/or natural variants of Galton-Watson
trees (see e.g. Jagers \cite{Jag-68}):
\begin{itemize}\item ${\rm GW}(q)$-trees as the most basic model just described;
  \item ${\rm GW}(q,\kappa)$-trees as ${\rm GW}(q)$-trees, where each individual is marked by an independent
    \em lifetime \em with distribution $\kappa$ on $(0,\infty)$; this includes the case
    $\kappa={\rm Exp}(c)$ of the exponential lifetime distribution with rate parameter
    $c\in(0,\infty)$;
  \item ${\rm GW}(q,\kappa,\beta)$-bushes as \em bushes \em (finite sequences) of a random number $N$
    of ${\rm GW}(q,\kappa)$-trees, where $N$ is Poisson distributed with parameter $\beta\in[0,\infty)$, in shorthand: $N\sim{\rm Poi}(\beta)$;
  \item ${\rm GWI}(q,\kappa,\eta,\chi)$-forests as \em forests \em (point processes on the \em forest
    floor \em $[0,\infty)$) of independent bushes of $N_i$ ${\rm GW}(q,\kappa)$-trees at the locations
    $S_i$ of a renewal process with inter-renewal distribution $\chi$ on $(0,\infty)$, where
    each $N_i$, $i\ge 1$, has distribution $\eta$ on $\bN=\{1,2,\ldots\}$.
\end{itemize}
With the common interpretation that individuals give birth only at the end of their life and that renewal locations are \em immigration times\em, all but the first model give rise to continuous-time processes counting the number $Y_t$ of individuals in the population at
time $t\ge 0$. In these continuous-time models it is natural to take $q(1)=0$, since an
individual producing a single child at its death time can be viewed as continuing to live instead of being replaced by its child.

\begin{figure}[ht]
\psfrag{Tsub}{$T_{\rm sub}^{p-\rm rdc}$}
\psfrag{Tcol}{$T^{p-\rm col}$}
\psfrag{Tblue}{$T^{p-\rm rdc}$}
\psfrag{time}{\scalebox{.7}{$t$}}

\epsfxsize=7cm
\centerline{\epsfbox{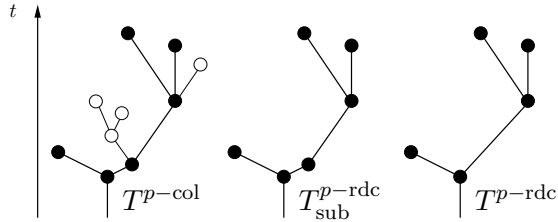}}
\caption{{Black vertices are represented by solid circles
and red ones by circle lines. 
%The dashed arrows represent a
%   (re)construction procedure to transform $\tau_1\sim\tau^{p\rm-rdc}$ (back) into $\tau_3\sim\tau^{\rm col}$.
}}
\label{discperc}
\end{figure}

Reduction by Bernoulli leaf colouring was studied in \cite{DuW} and reads as follows in our
setting:
%\begin{minipage}[t]{9cm}
\begin{itemize}\item independently mark each leaf of a tree $T$ (with lifetimes), or of a bush $B=(T_{(1)},\ldots,T_{(N)})$ or of a forest $F=(B(t),t\ge 0)$ by a
  \em black \em colour mark with probability $1-p\in(0,1)$, by a \em red \em colour mark otherwise; for this to be non-trivial, let $q(0)>0$, as $T$ then has leaves;
  \item if there are any black leaves, consider, as illustrated in Figure \ref{discperc},
  \begin{itemize}\item the \em $p$-reduced subtree \em $T^{p-\rm rdc}_{\rm sub}$ of $T$ as the
    subtree of $T$ spanned by the root and the black leaves (with lifetime marks inherited);
  \item and the \em $p$-reduced tree \em $T^{p-\rm rdc}$ derived from $T^{p-\rm rdc}_{\rm sub}$ by
    identifying vertices via the equivalence relation generated by $v\equiv w$ for vertices
    in $T_{\rm sub}^{p-\rm rdc}$ if $v\rightarrow w$ and $w$ is the only child of $v$ in $T_{\rm sub}^{p-\rm rdc}$
    (marked by the sum of lifetimes in each equivalence class);
  \item or the \em $p$-reduced bush \em
    $B^{p-\rm rdc}=(T_{(I_1)}^{p-\rm rdc},\ldots,T_{(I_{N^{p-\rm rdc}})}^{p-\rm rdc})$ as the
    $p$-reduced trees associated with the subsequence $(I_1,\ldots,I_{N^{p-\rm rdc}})$ of trees
    in $B$ that have black leaves;
  \item or the \em $p$-reduced forest \em $F^{p-\rm rdc}=(B^{p-\rm rdc}(t),t\ge 0)$ of
    $p$-reduced bushes.
  \end{itemize}
\end{itemize}
\vspace{-0.3cm}

%\end{minipage}\hfill
%\begin{minipage}[t]{7cm}

%\begin{figure}[ht]
%\psfrag{taugw}{$\tau^{\rm col}$}
%\psfrag{tausub}{$\tau^{p\rm-rdc}_{\rm sub}$}
%\psfrag{taublack}{$\tau^{p\rm-rdc}$}
%\psfrag{tau1}{$\tau_1 $}
%\psfrag{tau2}{$\tau_2 $}
%\psfrag{tau3}{$\tau_3 $}
%
%\epsfxsize=5cm
%\centerline{\epsfbox{figure1b3.eps}}
%
%\end{figure}
%\end{minipage}
It is easily seen that if $T$ is a Galton-Watson tree, then given that there are
any black leaves, the $p$-reduced subtree $T_{\rm sub}^{p-\rm rdc}$ and the $p$-reduced tree $T^{p-\rm rdc}$ are also Galton-Watson trees \cite{DuW}. By standard thinning properties of Poisson point processes, the $p$-reduced bush $B^{p-\rm rdc}$ associated with a ${\rm GW}(q,\kappa,\beta)$-bush $B$ is also a
Galton-Watson bush. We refer to the offspring distribution $q^{p-\rm rdc}$, the lifetime
distribution $\kappa^{p-\rm rdc}$ and the Poisson parameter $\beta^{p-\rm rdc}$ of $B^{p-\rm rdc}$ as the \em
$p$-reduced triplet $(q^{p-\rm rdc},\kappa^{p-\rm rdc},\beta^{p-\rm rdc})$ associated with $(q,\kappa,\beta)$\em, similarly for forests etc. It is not hard to find offspring distributions $\widehat{q}$ that do \em not \em arise as $p$-reduced offspring distributions for any $q$ (e.g. $\widehat{q}(0)=\widehat{q}(3)=1/2$). More precisely, \cite{DuW} obtained the following characterisation.
\begin{theo}[Theorem 4.2 of \cite{DuW}]\label{thm1} For an offspring distribution $q$, the following are
  equivalent:
  \begin{enumerate}\item[\rm(i)] There is a family $(q_\lambda)_{\lambda\ge 0}$ of offspring
    distributions with $q_1=q$ such that $q_\mu$ is the $(1-\mu/\lambda)$-reduced offspring
    distribution associated with $q_\lambda$, for all $0\le\mu<\lambda<\infty$.
    \item[\rm(ii)] The generating function $\varphi_q$ of $q$ satisfies
      \begin{equation} \varphi_q(s)=\sum_{i=0}^\infty s^iq(i)=s+\widetilde{\psi}(1-s),\qquad 0\le s\le 1,\label{genfnbm}
      \end{equation}
      where for some $\widetilde{b}\in\bR$, $\widetilde{a}\ge 0$ and a measure $\widetilde{\Pi}$ on $(0,\infty)$ with $\int_{(0,\infty)}(1\wedge x^2)\widetilde{\Pi}(dx)<\infty$,
      \begin{equation} \widetilde{\psi}(r)=\widetilde{b}r+\widetilde{a}r^2+\int_{(0,\infty)}(e^{-rx}-1+rx\mathbf{1}_{\{x<1\}})\widetilde{\Pi}(dx),\qquad r\ge 0.\label{lkbm}
      \end{equation}
  \end{enumerate}
  In the setting of {\rm (i)} and {\rm (ii)}, a consistent family
  $(B_\lambda)_{\lambda\ge 0}$ of
  ${\rm GW}(q_\lambda,{\rm Exp}(c_\lambda),\beta_\lambda)$-bushes can be constructed such that $(B_\mu,B_\lambda)\overset{{\rm (d)}}{=}(B_\lambda^{(1-\mu/\lambda)-\rm rdc},B_\lambda)$ for all
  $0\le\mu<\lambda<\infty$. For each $c=c_1\in(0,\infty)$ and $\beta=\beta_1\in(0,\infty)$, the family $(q_\lambda,c_\lambda,\beta_\lambda)_{\lambda\ge 0}$ is unique, $(q_\lambda)_{\lambda\ge 0}$ does not depend on $(c,\beta)$, while $(c_\lambda)_{\lambda\ge 0}$ depends on $q$ but not on $\beta$ and $(\beta_\lambda)_{\lambda\ge 0}$ on $q$ but not on $c$.
\end{theo}
In \cite{DuW}, this result is a key step in the construction of L\'evy trees as
genealogies of Markovian continuous-state branching processes with branching mechanism $\psi$, where $\psi$ is a linear transformation of $\widetilde{\psi}$ that we recall in Section \ref{Consistent family}. In the present paper we establish characterisations analogous to Theorem \ref{thm1} for the other variants of
Galton-Watson trees, bushes and forests.

\begin{theo}\label{thm2} For a pair $(q,\kappa)$ of offspring and lifetime distributions, the
  following are equivalent:
  \begin{enumerate}\item[\rm(i)] There are families $(q_\lambda,\kappa_\lambda)_{\lambda\ge 0}$
    with $q_1=q$ and $\kappa_1=\kappa$ such that
    $(q_\mu,\kappa_\mu)$ is the $(1-\mu/\lambda)$-reduced pair associated with
    $(q_\lambda,\kappa_\lambda)$, for all $0\le\mu<\lambda<\infty$.
    \item[\rm(ii)] The generating function $\varphi_q$ of $q$ satisfies $\varphi_q(s)=s+\widetilde{\psi}(1-s)$,
      where $\widetilde{\psi}$ is of the form {\rm(\ref{lkbm})}. Moreover, $\kappa$ is geometrically
      divisible in that there is a family $(X_\alpha^{(j)},j\ge 1)$ of independent identically
      distributed random variables and $G^{(\alpha)}\sim{\rm Geo}(\alpha)$ independent geometric with parameter
      $\alpha$, i.e. $\bP(G^{(\alpha)}=k)=\alpha(1-\alpha)^{k-1}$, $k\in\bN$, such that $X_\alpha^{(1)}+\cdots+X_\alpha^{(G^{(\alpha)})}\sim\kappa$
      \begin{itemize}\item for all $\alpha>1/\widetilde{\psi}^\prime(\infty)$ if $\widetilde{\psi}^\prime(\infty)<\infty$, where $\widetilde{\psi}^\prime(\infty)$ means $\lim_{r\rightarrow\infty}\widetilde{\psi}^\prime(r)$;
        \item for all $\alpha>0$ if $\widetilde{\psi}^\prime(\infty)=\infty$.
      \end{itemize}
  \end{enumerate}
  In the setting of {\rm (i)} and {\rm (ii)}, a consistent family $(B_\lambda)_{\lambda\ge 0}$ of
  ${\rm GW}(q_\lambda,\kappa_\lambda,\beta_\lambda)$-bushes can be constructed such that $(B_\mu,B_\lambda)\overset{{\rm (d)}}{=} (B_\lambda^{(1-\mu/\lambda)-\rm rdc},B_\lambda)$ for all
  $0\le\mu<\lambda<\infty$. For each $\beta=\beta_1\in(0,\infty)$, the family $(q_\lambda,\kappa_\lambda,\beta_\lambda)_{\lambda\ge 0}$ is unique, $(q_\lambda)_{\lambda\ge 0}$ does not depend on $(\kappa,\beta)$ while $(\kappa_\lambda)_{\lambda\ge 0}$ depends on $q$ but not on $\beta$ and $(\beta_\lambda)_{\lambda\ge 0}$ on $q$ but not on $\kappa$.
\end{theo}
The requirement on $\kappa$ set in the second bullet point is referred to as geometric infinite
divisibility in the literature, see \cite{KMM}, also Section \ref{contstate} here. Since the distribution $\kappa={\rm Exp}(c)$ is geometrically infinitely
divisible, Theorem \ref{thm2} is an extension of Theorem \ref{thm1}.

\begin{theo}\label{thm3} For a pair $(q,\eta)$ of offspring and immigration distributions, the following are equivalent:
  \begin{enumerate}\item[\rm(i)] There are families $(q_\lambda,\eta_\lambda)_{\lambda\ge 0}$
    with $q_1=q$ and $\eta_1=\eta$ such that
    $(q_\mu,\eta_\mu)$ is the $(1-\mu/\lambda)$-reduced pair associated with
    $(q_\lambda,\eta_\lambda)$, for all $0\le\mu<\lambda<\infty$.
    \item[\rm(ii)] The generating function $\varphi_q$ of $q$ satisfies $\varphi_q(s)=s+\widetilde{\psi}(1-s)$,
      where $\widetilde{\psi}$ is of the form {\rm(\ref{lkbm})}. Moreover, the generating function $\varphi_\eta$ of $\eta$ satisfies
    $$\varphi_\eta(s)=\sum_{i=1}^\infty s^i\eta(i)=1-\widetilde{\phi}(1-s),\qquad 0\le s\le 1,\label{genfnim}
    $$
      where for some $\widetilde{d}\in\bR$, and a measure $\widetilde{\Lambda}$ on $(0,\infty)$ with $\int_{(0,\infty)}(1\wedge x)\widetilde{\Lambda}(dx)<\infty$,
      \begin{equation} \widetilde{\phi}(r)=\widetilde{d}r+\int_{(0,\infty)}(1-e^{-rx})\widetilde{\Lambda}(dx),\qquad r\ge 0.\label{lkim}
      \end{equation}
  \end{enumerate}
  In the setting of {\rm (i)} and {\rm (ii)}, a consistent family $(F_\lambda)_{\lambda\ge 0}$ of
  ${\rm GWI}(q_\lambda,{\rm Exp}(c_\lambda),\eta_\lambda,{\rm Exp}(\paramgamma_\lambda))$-forests can
  be constructed such that $(F_\mu,F_\lambda)\overset{{\rm (d)}}{=} (F_\lambda^{(1-\mu/\lambda)-\rm rdc},F_\lambda)$ for all
  $0\le\mu<\lambda<\infty$. For each $c=c_1\in(0,\infty)$ and $\paramgamma=\paramgamma_1\in(0,\infty)$, the family $(q_\lambda,c_\lambda,\eta_\lambda,\paramgamma_\lambda)_{\lambda\ge 0}$ is unique, $(q_\lambda)_{\lambda\ge 0}$ does not depend on $(c,\eta,\paramgamma)$, while $(c_\lambda)_{\lambda\ge 0}$ depends on $q$ but not on $(\eta,\paramgamma)$, $(\eta_\lambda)_{\lambda\ge 0}$ depends on $q$ but not on $(c,\paramgamma)$ and $(\paramgamma_\lambda)_{\lambda\ge 0}$ depends on $(q,\eta)$ but not on $c$. 
\end{theo}
\noindent The binary special case with single immigrants, where for some $\theta\ge 0$ and all $\lambda\ge 0$ 
$$q_\lambda(0)=\frac{1}{2}+\frac{1}{2\sqrt{\theta^2+2\lambda}},\quad q_\lambda(2)=\frac{1}{2}-\frac{1}{2\sqrt{\theta^2+2\lambda}},\qquad c_\lambda=\sqrt{\theta^2+2\lambda},$$
$$\eta_\lambda(1)=1,\qquad h_\lambda=\sqrt{\theta^2+2\lambda}-\theta,$$
leads to the setting of \cite{PiW-05}, where $(F_\lambda)_{\lambda\ge 0}$ was shown to have independent ``increments'' expressed by a composition rule, and to converge to the forest in Brownian motion with drift $-\theta$. 

Theorems \ref{thm2} and \ref{thm3} describe in the same way the genealogy of associated continuous-state branching processes (CSBP) as Theorem \ref{thm1}. Specifically, for Theorem
\ref{thm2} the continuous-state processes are Sagitov's age-dependent ${\rm CSBP}(K,\psi)$ based on a branching mechanism $\psi$ and the distribution of a local time process $K$, i.e.\ either an inverse subordinator or an inverse increasing random walk, see \cite{KaS-98,Sag-91} and Section \ref{CSBPKpsi} here; for Theorem \ref{thm3}, they are CSBP with immigration,
${\rm CBI}(\psi,\phi)$, where $\phi$ is an immigration mechanism, see \cite{KaW-71,Lam-CBI} and Section \ref{contstate} here.

\begin{prop}\label{prop4a} Let $(Z_t^\lambda,t\ge 0)$ be the population size process in the setting of Theorem \ref{thm2}. If $\widetilde{\psi}^\prime(0)>-\infty$, then
$$ \frac{Z^\lambda_t}{\psi^{-1}(\lambda)}\rightarrow Z_t\qquad\mbox{almost surely as $\lambda\rightarrow\infty$, for all $t\ge 0$,}
$$
where $(Z_t,t\ge 0)$ is a ${\rm CSBP}(K,\psi)$ with $Z_0=\beta$, for some $\psi$ linear transformation of $\widetilde{\psi}$ and $K=(K_s,s\ge 0)$ such that $\inf\{s\ge 0\colon K_s>V_\lambda\}\sim\kappa_\lambda$ for $V_\lambda\sim{\rm Exp}(c_\lambda)$ with $c_\lambda$ as in Theorem \ref{thm1}. 
%If furthermore $\psi^\prime(0)>-\infty$, then the convergence holds in the almost sure sense.
\end{prop}

\begin{prop}\label{prop4} Let $(Y_t^\lambda,t\ge 0)$ be the population size process in the setting of Theorem \ref{thm3}. Then
$$ \frac{Y_t^\lambda}{\psi^{-1}(\lambda)}\rightarrow Y_t\qquad\mbox{in distribution as $\lambda\rightarrow\infty$, for all $t\ge 0$,}
$$
where $(Y_t,t\ge 0)$ is a ${\rm CBI}(\psi,\phi)$ with $Y_0\!=\!0$, for $\psi$ and $\phi$ linear transformations of $\widetilde{\psi}$ and $\widetilde{\phi}$. If furthermore $\psi^\prime(0)>-\infty$ and $\phi^\prime(0)<\infty$, then the convergence holds in the almost sure sense.
\end{prop}

\noindent These convergence results should be seen in the context of the large literature on space-time scaling limits of branching processes in discrete or continuous time, see
\cite{Duq-imm,DuW2,KaS-98,KaW-71,Lam-lim,Pin-lim,Sag-91}. Convergence in distribution holds under much weaker assumptions on the families $(q_\lambda,\kappa_\lambda,\beta_\lambda)_{\lambda\ge 0}$ or
$(q_\lambda,c_\lambda,\eta_\lambda,h_\lambda)_{\lambda\ge 0}$ and invariance principles in varying degrees of generality have been obtained. It is also well-known that the convergence in distribution at a fixed time for Markovian branching processes implies the convergence in distribution of the whole process in the Skorohod sense of convergence of right-continuous functions with left limits. In \cite{DuW2}, joint convergence of processes and their genealogical trees is shown, also for a wider class of families $(B_\lambda)_{\lambda\ge 0}$ suitably converging to bushes of L\'evy trees. The main  contribution of the present work is to provide almost sure approximations of more general classes of continuous-state processes and consistent families of trees that contain full information about the genealogy of the population of the limiting continuous-state process, which is not contained in the limiting process itself nor in the approximating discrete-state branching processes.

The structure of this paper is as follows. In Section \ref{preliminaries} we formally set up the framework in which we represent trees, we recall preliminaries from Duquesne and Winkel \cite{DuW} and develop a bit further some aspects that readily transfer and serve in the more general context here. We also provide some background about continuous-state branching processes with immigration, and about geometric infinite divisibility. Section \ref{Consistend kappa bush} presents the theory around Theorem \ref{thm2} and Proposition \ref{prop4a}, while Section \ref{Section Immigration} deals with Theorem \ref{thm3} and Proposition \ref{prop4}. In each setting we provide explicit
formulas for offspring distributions, lifetime distributions and immigration distributions
as appropriate; we also provide explicit reconstruction procedures that reverse the reduction for the consistent families
of bushes and forests and establish connections with backbone decompositions (Theorem \ref{CSBPbackbone}) and L\'evy trees. We finally deduce generalisations combining lifetimes and immigration.

%\section{Preliminaries}\label{prel}

\section{Preliminaries}\label{preliminaries}
%This section introduces the basic notions and operations on the
%trees that we use. We follow Neveu \cite{jn:nv tree} ( see also
%Chauvin \cite{cha:pro} ). In this tree model we will recall from
%Duquesne and Winkel \cite{mw&td:lv tree} the reductions and growth
%of GW trees. We also introduce continuous-time branching processes
%with discrete and continuous state spaces. Some more definitions
%such as subordinators, geometric infinite divisibility,
%continuous-state branching process (CSBP) and CSBP with immigration
%(CBI) process will be given.
%--------------------------------------------------------------------------------------------------------------------------------------Discrete trees and GW branching property
\subsection{Discrete trees with edge lengths and colour marks}\label{trees}
\subsubsection{Discrete trees and the Galton-Watson branching
property}\label{Discr Trees and GW BP}
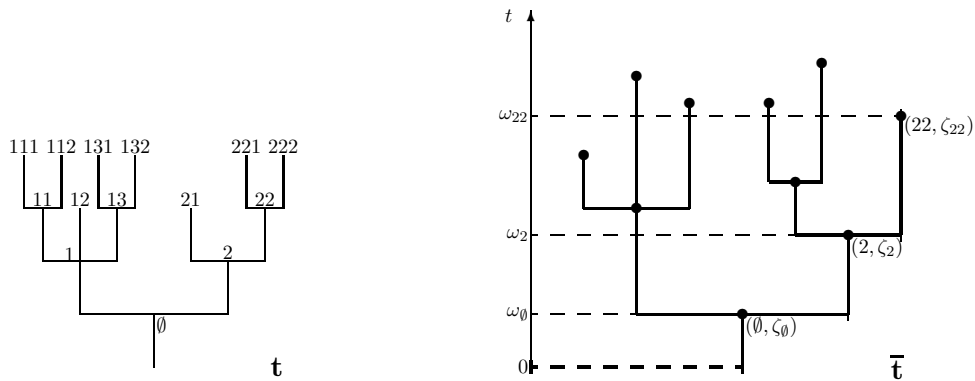
\begin{figure}[b]
%\begin{minipage}{3cm}
\begin{center}
\begin{picture}(100,75)(0,-10)
\put(50,-10){\line(0,1){20}} \put(22,10){\line(50,0){56}}
\put(22,10){\line(0,1){20}}
\put(8,30){\line(0,1){20}}\put(8,30){\line(1,0){28}}
\put(1,50){\line(1,0){14}} \put(15,50){\line(0,1){20}}
\put(1,50){\line(0,1){20}}\put(36,30){\line(0,1){20}}
\put(29,50){\line(1,0){7}}\put(29,50){\line(0,1){20}}

\put(78,10){\line(0,1){20}}\put(78,30){\line(1,0){14}}\put(92,30){\line(0,1){20}}
\put(92,50){\line(1,0){7}}\put(99,50){\line(0,1){20}}
\put(22,30){\line(0,1){20}}\put(36,50){\line(1,0){7}}\put(43,50){\line(0,1){20}}
\put(78,30){\line(-1,0){14}}\put(64,30){\line(0,1){20}}\put(92,50){\line(-1,0){7}}
\put(85,50){\line(0,1){20}}

\put(20,35){\makebox(0,0)[tr]{\scalebox{.7}{1}}}
\put(15,71){\makebox(0,0)[b]{\scalebox{.7}{112}}}
\put(8,51){\makebox(0,0)[b]{\scalebox{.7}{11}}}
\put(78,31){\makebox(0,0)[b]{\scalebox{.7}{2}}}
\put(1,71){\makebox(0,0)[b]{\scalebox{.7}{111}}}
\put(36,51){\makebox(0,0)[b]{\scalebox{.7}{13}}}
\put(92,51){\makebox(0,0)[b]{\scalebox{.7}{22}}}
\put(99,71){\makebox(0,0)[b]{\scalebox{.7}{222}}}
\put(22,51){\makebox(0,0)[b]{\scalebox{.7}{12}}}
\put(43,71){\makebox(0,0)[b]{\scalebox{.7}{132}}}
\put(64,51){\makebox(0,0)[b]{\scalebox{.7}{21}}}
\put(85,71){\makebox(0,0)[b]{\scalebox{.7}{221}}}
\put(29,71){\makebox(0,0)[b]{\scalebox{.7}{131}}}
\put(51,9){\makebox(0,0)[tl]{\scalebox{.7}{$\emptyset$}}}
\put(99,-10){\makebox(0,0)[r]{$\mathbf{t}$}}
\end{picture}
\hspace{3cm}
%\begin{minipage}{6cm}
\begin{picture}(160,140)(0,-20)
\put(0,0){\dashbox{5}(120,0)}
\put(0,30){\dashbox{5}(140,0)}\put(0,75){\dashbox{5}(140,0)}
\put(0,-20){\vector(0,1){135}} \thicklines
\put(80,-20){\line(0,1){20}}\put(40,0){\line(1,0){80}}
\put(0,-20){\dashbox{5}(80,0)}
\put(40,0){\line(0,1){90}}\put(20,40){\line(1,0){40}}\put(60,40){\line(0,1){40}}
\put(20,40){\line(0,1){20}}\put(120,0){\line(0,1){30}}\put(100,30){\line(1,0){40}}
\put(100,30){\line(0,1){20}}\put(90,50){\line(1,0){20}}\put(90,50){\line(0,1){30}}
\put(110,50){\line(0,1){45}}\put(140,30){\line(0,1){45}}

\put(-7,113){\makebox(0,0)[r]{\scalebox{.7}{$t$}}}\put(120,30){\circle*{4}}
\put(141,-20){\makebox(0,0)[r]{$\overline{\mathbf{t}}$}}
\put(121,29){\makebox(0,0)[tl]{\scalebox{.7}{$(2,\zeta_2)$}}}
\put(141,75){\makebox(0,0)[tl]{\scalebox{.7}{$(22,\zeta_{22})$}}}
\put(-1,0){\makebox(0,0)[r]{\scalebox{.7}{$\omega_\emptyset$}}}
\put(-1,-20){\makebox(0,0)[r]{\scalebox{.7}{0}}}
\put(-1,30){\makebox(0,0)[r]{\scalebox{.7}{$\omega_2$}}}
\put(-1,75){\makebox(0,0)[r]{\scalebox{.7}{$\omega_{22}$}}}
\put(81,-1){\makebox(0,0)[tl]{\scalebox{.7}{$(\emptyset,\zeta_\emptyset)$}}}
\put(40,90){\circle*{4}}\put(20,60){\circle*{4}}\put(80,0){\circle*{4}}\put(60,80){\circle*{4}}
\put(100,50){\circle*{4}}%\put(100,85){\circle*{4}}
\put(40,40){\circle*{4}}\put(90,80){\circle*{4}}
\put(110,95){\circle*{4}}\put(140,75){\circle*{4}}

\end{picture}
\end{center}
%\end{minipage}
\vspace{-0.3cm}
\caption{{On the left, $\ft=\{\emptyset, 1,2,11,12,13,21,22,111,112,131,132,221,222\}$, and on the right
$\overline{\ft}=\{(\emptyset,\zeta_\emptyset),(1,\zeta_1),(11,\zeta_{11}),(12,\zeta_{12}),(13,\zeta_{13}),%$\newline %\hfill$\;$
(2,\zeta_2),(21,\zeta_{21}),(211,\zeta_{211}),(212,\zeta_{212}),(22,\zeta_{22})\}$.}}
%D_2=t_1,$ $D_{22}=t_2,$ $x_2=t_1-t_0,$ $x_{22}=t_2-t_1.D_\emptyset = x_\emptyset = t_0$}

%\end{minipage}
\label{fig3}
\label{fig2}
\end{figure}

Following Neveu \cite{Nev-GW}, Chauvin \cite{Cha-BH} and others, we let
$$
\bU=\bigcup_{n \geq 0}
\bN^n=\{\emptyset,1,2,\ldots,11,12,\ldots,21,22,\ldots,111,112,\ldots\}
$$
be the set of \em integer words\em, where $\mathbb{N}=\{1,2,3,\ldots \}$ and
where $\mathbb{N}^{0}=\{\emptyset\}$ has the \em empty word \em $\emptyset$ as its only element.
For $u=u_1u_2\cdots u_n\in\bU,$ and $v=v_1v_2\cdots v_m\in\bU$, we denote by
$|u|$ the \em length \em of $u$, e.g. $|u_1u_2\cdots u_n|=n$, and by
$uv=u_1u_2\cdots u_nv_1v_2\cdots v_m$ the \em concatenation \em of words
in $\bU$. %We define a tree $\mathbf{t}$ as a subset of
%$\mathbb{U}$ with certain properties as follows:
\begin{defi}\label{defn1}
\rm A subset $\mathbf{t}\subset\bU$ is called a \em tree \em if
\begin{itemize}
    \item $\emptyset \in \mathbf{t}$; we refer to $\emptyset$ as the \em progenitor \em of $\ft$;
    \item for all $u\in \mathbb{U}$ and $j \in \mathbb{N}$ with $uj\in \mathbf{t}$, we have $u \in \mathbf{t}$; we refer to $u$ as the \em parent \em of $uj$;
    \item for all $u\in \mathbf{t}$, there exists $\nu_{u}(\mathbf{t})\in\bN_0=\{0,1,2,\ldots\}$ such that $uj\in \mathbf{t} \Leftrightarrow 1\leq j \leq
    \nu_{u}(\mathbf{t})$; we refer to $\nu_{u}(\mathbf{t})$ as the \em number of children \em of $u$.
\end{itemize}
\end{defi}
We refer to the length $|u|$ of a word $u\in\ft$ as the \em generation \em of the \em individual \em $u$ in the \em genealogical tree \em $\ft$. An element $u\in \mathbf{t}$ is called a
\em leaf \em of $\mathbf{t}$ if and only if $\nu_{u}(\mathbf{t})=0$. We consider the lexicographical order $\le$ on $\bU$ and its restriction to $\mathbf{t}$ as the canonical \em total order\em. For $u,v\in\bU$, we write $u
\preceq uv$, defining a \em partial order \em on $\bU$, whose restriction to $\ft$ is the \em genealogical order \em on $\ft$. The partial order $\preceq$ is compatible with the total order $\le$ in that $u\preceq v\Rightarrow u\le q$. A tree $\ft$ in the sense of Definition \ref{defn1} can be represented graphically as in Figure \ref{fig2}.

%\begin{center}
%\end{center}
Let $\mathbb{T}$ be the space of all such trees, and let
$\mathbb{T}_{u}=\{\mathbf{t} \in \mathbb{T}\colon 
u \in \mathbf{t} \}$. Then $\nu_u$ is a map defined on $\bT_u$ taking values in
$\bN_0$. Note that $\mathbb{T}$ is
uncountable. A sigma-algebra on $\mathbb{T}$ is defined as $\mathcal{F}=\sigma\{\mathbb{T}_{u},u\in\bU\}$. We also specify the $n$th generation
$\pi_n(\ft)=\{u\in\ft\colon |u|=n\}=\ft\cap\bN^n$ and set
$\mathcal{F}_{n}=\sigma\{\mathbb{T}_{u}, |u|\leq n\}=\sigma\{\pi_m,m\le n\}$.%, set $\bT_A=\bigcap_{u\in A}\bT_u$ for finite $A\subset\bN^n$.

We define the shift map/operator $\theta_u$ that assigns to a tree $\ft$ its \em subtree $\ft_u=\theta_u\ft$ above $u \in\ft$\em:
$$\theta_u\colon  \mathbb{T}_{u} \rightarrow
\mathbb{T}, \qquad \mathbf{t} \mapsto \theta_u\ft=\ft_u=\{v \in
\mathbb{U}\colon  uv \in \mathbf{t}\}.$$
Clearly $\bT_\emptyset=\bT$ and
$\{\nu_u \geq j\}\cap \mathbb{T}_{u}= \mathbb{T}_{uj}$, also
$\theta_u^{-1}(\mathbb{T}_{v})=\mathbb{T}_{uv}$ and
%and for $v=v_1v_2\cdots v_m$,
$$\mathbb{T}_{v}=\{\mathbf{t} \in \mathbb{T}\colon \nu_{v_1v_2\cdots v_k}(\mathbf{t})\geq
{v_{k+1}}\mbox{ for all $0\leq k < m$}\}\qquad\mbox{for $v=v_1v_2\cdots v_m\in\bU$.}$$
These relations allow us to consider a \em random \em tree $\tau$ whose distribution is a
probability measure $\bP_q$ on $\mathbb{T}$, under which the
numbers of children $\nu_u$ of the individuals $u$ in the random tree are independent random variables with distribution $q$. More formally, $\bP_q$ is
characterised as follows:
%The Galton-Watson branching property states that the  shifted trees
%starting from $u \in \mathbf{t}, |u|=n,$ are i.i.d. copies of the
%original tree
%$\mathbf{t}$.\\
%\newpage
\paragraph{${\rm GW}(q)$-trees and their branching property (see e.g. Neveu \cite{Nev-GW})}
\begin{itemize}
    \item[(a)] For any probability measure $q$ on $\mathbb{N}_0$, there exists a unique
      probability measure $\mathbb{P}_q$ on $(\mathbb{T}, \mathcal{F})$ such that
      $\mathbb{P}_q(\nu_{\emptyset}=j)=q(j)$, and conditionally on $\{\nu_{\emptyset}=j\}$ for
      any $j\ge 1$ with $q(j)>0$, the subtrees $\theta_{i}$, $1\le i\leq j$, above the
      first generation are independent with distribution $\mathbb{P}_q$. A random tree
      $\tau$ with distribution $\mathbb{P}_q$ is called a ${\rm GW}(q)$-tree.
    \item[(b)] Under $\bP_q(\,\cdot\,|\cF_n,\pi_n=A)$,
      the subtrees $\theta_u$, $u\in A$, above the $n$th generation are
      independent and with distribution $\mathbb{P}_q$, for all finite $A\subset\bN^n$ and $n\ge 1$ with $\bP(\pi_n=A)>0$.
\end{itemize}
For finite trees, in particular in the \em (sub)critical \em case $\bE_q(\nu_\emptyset)=\sum_{j\in\bN_0}jq(j)\le 1$, the measure $\bP_q$ can be expressed as
$$ \mathbb{P}_q(\{\mathbf{t}\})=\prod_{v \in \mathbf{t}} q(\nu_v(\mathbf{t})), \qquad \mbox{for all $\mathbf{t} \in {\mathbb{T}}$,}$$
but this does not specify the measure $\bP_q$ in the \em supercritical \em case  $\bE_q(\nu_\emptyset)>1$, where $\bP_q$ assigns positive measure to infinite trees. Here, $\bE_q$ is the expectation operator associated with $\bP_q$.
For a ${\rm GW}(q)$-tree $\tau$, the process $G_n=\#\pi_n(\tau)$, $n\ge 0$, is known as a ${\rm GW}(q)$-branching process.

%-----------------------------------------------------------------------------------------------------------------Lifetime marks and discrete BP in cont. times
\subsubsection{Marked trees and discrete branching processes in continuous time}\label{Lifetime marks and discr BP in cont. time}
Let $(\bH,\cH)$ be a measurable space of marks. We can attach a mark $\xi_u\in\bH$ to each vertex $u$ of a given tree $\mathbf{t}$. Formally, a marked tree is a subset
\begin{equation}\overline{\mathbf{t}}\subset \mathbb{U}\times \mathbb{H} \quad
\mbox{such that} \quad \mathbf{t}= \{u\in\bU\colon (u,\xi_u)\in\overline{\mathbf{t}}\mbox{ for some $\xi_u\in\bH$}\}\in
\mathbb{T}\label{skel}\end{equation}
%----
and where $\overline{\ft}\cap\{u\}\times\bH=\{(u,\xi_u)\}$, i.e. the map $\xi\colon \ft\rightarrow\bH$ is unique. So a marked tree has the form $\overline{\ft}=\{(u,\xi_u)\in\bU\times\bH\colon u\in \mathbf{t}\}$. We denote by
$\bT^\bH$ the set of marked trees. We consider the set of trees $\mathbb{T}_{u}^\bH=\{\overline{\mathbf{t}} \in\bT^\bH\colon  u \in \mathbf{t}\}$ containing
individual $u$ and note that $\xi_u\colon \bT_u^\bH\rightarrow\bH$ is a map.

For marked trees, we set $\nu_u(\overline{\mathbf{t}})=\nu_u(\mathbf{t})$ and
$\overline{\ft}_u=\overline{\theta}_u\overline{\ft}=\{(v,\xi_v)\colon v
\in \mathbf{t}_u\}=\{(v,\xi_{v})\colon uv \in \mathbf{t}\}$. As sigma-algebra on $\bT^\bH$ we take
one that makes $\overline{\ft}\mapsto\ft$ in (\ref{skel}) and $\overline{\ft}\mapsto\xi_u(\overline{\ft})$ measurable:
$$\mathcal{F}^{\bH}
%=\sigma\{\overline{\mathbb{T}}_{u}, u\in \mathbb{U}\}\vee \sigma\{\xi_u, u\in \mathbb{U}\}
=\sigma\{\bT_{u,H}^\bH,u\in\bU,H\in\cH\},\mbox{ where $\bT_{u,H}^\bH=\{\overline{\ft}\in\mathbb{T}_u^\bH\colon \xi_u\in H\}$.}$$

For example, for $\bH=(0,\infty)$, the marks can represent the \em lifetimes \em of individuals. We will later use $\bH=(0,\infty)\times\{0,1\}$ so that $\xi_u=(\zeta_u,\gamma_u)$
consists of a lifetime mark $\zeta_u\in(0,\infty)$ and a colour mark $\gamma_u\in\{0,1\}$.
%In \cite{CeG,CGM}, $u$ is marked by a lifetime $\zeta_u \in[0, \infty]$ with the assumption that if
%$\zeta_u=\infty,$ then $\nu_u=0$. 
We consider a model when
$u \in \mathbf{t}$ will produce children at the moment of its death.
The birth and death times $\alpha_u$ and $\omega_u$ of each individual $u \in \mathbf{t},$ are then defined recursively by
$$\left\{ \begin{array}{ll}
         \alpha_\emptyset=0,\ \omega_{\emptyset}=\zeta_\emptyset,\\
         \alpha_{uj}=\omega_u,\ \omega_{uj}=\alpha_{uj}+\zeta_{uj},&1\le j\le\nu_u, u\in\ft.
         \end{array}
  \right. $$
We denote by $\overline{\pi}_t(\overline{\ft})=\{u\in\mathbf{t}\colon \alpha_u<t\le\omega_u\}$ the set of individuals
alive at time $t\ge 0$ and define $\cF_t^\bH=\sigma\{\overline{\pi}_s,s\le t\}$. For $u\in\overline{\pi}_t(\overline{\ft})$, we denote by
$$ \overline{\theta}_{u,t}(\overline{\ft})=\{(\emptyset,\omega_u-t)\}\cup\{(v,\zeta_{uv})\colon uv\in\ft\}
$$
the subtree of individual $u$ above $t$. Figure \ref{fig3} gives a graphical representation of a tree $\overline{\ft}$ where lifetimes are shown as edge lengths.

%\begin{figure}[h]

%\end{figure}
%\newpage

In this formalism, we can define and study ${\rm GW}(q,\kappa)$-trees as ${\rm GW}(q)$ trees with
independent lifetimes distributed according to a measure $\kappa$ on $(0,\infty)$:
%Let us use
%${\rm GW}(q, \kappa)$ for a Galton-Watson tree with offspring distribution
%$q$ and i.i.d lifetimes with distribution $\kappa$.
%An important special case is
%$$\zeta_u \sim \textrm{Exp}(c), \hfill \textrm{for all} \ u \in \mathbf{t}.$$
%We will use the notation ${\rm GW}(q, \rm {Exp} (c))$ for this special
%case.

\paragraph{${\rm GW}(q,\kappa)$-trees and their branching property (see Neveu \cite{Nev-GW}, Chauvin \cite{Cha-BH})} \label{bran2}
\begin{itemize}
  \item[(a)] For any probability measure $\mathbb{Q}$ on $\mathbb{N}_0 \times \bH$,
    there exists a unique probability measure $\mathbb{P}_{\mathbb{Q}}$ on
    $(\mathbb{T}^\bH,\mathcal{F}^\bH)$, such that
    $(\nu_\emptyset,\xi_\emptyset)\sim {\mathbb{Q}}$ and conditionally on
    $\{\nu_\emptyset = j,\xi_\emptyset\in H\}$ for any $j\ge 1$, $H\in\cH$, with $\mathbb{Q}(\{j\}\times H)>0$, the subtrees
    $\overline{\theta}_i$, $1\le i \le j$, are independent
    with distribution $\mathbb{P}_{\mathbb{Q}}$. For $\bH=(0,\infty)$ and
    $\bQ=q\otimes\kappa$ a random tree $T$ with distribution $\bP_\bQ$ is called a
    ${\rm GW}(q,\kappa)$-tree.
  \item[(b)] Under
    $\bP_{q\otimes\kappa}(\,\cdot\,|\cF_t^\bH,\overline{\pi}_t=A)$, the subtrees
    $\overline{\theta}_{u,t}$, $u\in A$, above time $t$
    are independent and distributed like $\overline{\theta}_{\emptyset,s}$ under
    $\bP_{q\otimes\kappa}(\,\cdot\,|\zeta_\emptyset>s)$, where $s=t-\alpha_u$ is the
    ($\cF_t^\bH$-measurable) age of $u$ at time $t$, for all finite $A\subset\bU$ with
    $\bP(\overline{\pi}_t=A)>0$.
\end{itemize}
For further details including strong branching properties, we refer to \cite{CeG,CGM}.
For a ${\rm GW}(q,\kappa)$-tree $T$, the process $Z_t=\#\overline{\pi}_t(T)$, $t\ge 0$, is known as a Bellman-Harris branching process \cite{BeH2,BeH1}. The Markovian special case for $\kappa={\rm Exp}(c)$ is also
called a continuous-time Galton-Watson process.

%\end{document}

%%-------------------------------------------------------------------------------------------------coloured leaves and a two-colour branching property
\subsubsection{Coloured leaves, coloured trees and a two-colours branching
property}\label{coloured leaves and 2-colour BP}

On a suitable probability space $(\Omega,\mathcal{A},\mathbb{P})$, let $T\colon (\Omega,
\mathcal{A}, \mathbb{P})\!\rightarrow\!(\mathbb{T}^{(0,\infty)},
\mathcal{F}^{(0,\infty)}, \mathbb{P}_{q\otimes\kappa})$ be a ${\rm GW}(q,\kappa)$-tree.
%with offspring distribution
%$$\mathbb{P}_q (\nu_\emptyset=i)=q(i),\ i \geq 0.$$
We assume $q(0)>0$, i.e.\ $T$ has leaves, and $q(1)=0$, as an individual producing a
single child can be viewed as continuing to live instead of being replaced by its child.

Following \cite{DuW}, we independently mark the \em leaves \em $u$ of $T$ with one of two colours, say \em red\em, $\bP(\gamma_u(T)=1|\nu_u(T)=0)=p$, or \em black\em, $\bP(\gamma_u(T)=0|\nu_u(T)=0)=1-p$,
for some given $p\in(0,1)$. It will be convenient to also mark each non-leaf individual in black if the subtree above it has at least one black leaf, red otherwise. Such a marked tree $T^{p\rm-col}$ is a random element of $\bT^\bH$ for $\bH=(0,\infty)\times\{0,1\}$. We denote its distribution by $\bP_{q\otimes\kappa}^{p\rm-col}$. Note that it is \em not \em of the form
$\bP_\bQ$ introduced in the previous section, because marks for non-leaf individuals will not be independent. We set
\begin{equation}\label{gp}
g(p)=\mathbb{P}_{q\otimes\kappa}^{p\rm-col}(\gamma_{\emptyset}=1)=\mathbb{P}(\mbox{$T^{p\rm-col}$ has only red colour marks})
    =\mathbb{E}[p^{\#\{u\in T\colon \nu_u(T)=0\}}].
\end{equation}

\paragraph{Branching properties of coloured ${\rm GW}(q,\kappa)$-trees (cf. Duquesne and Winkel \cite{DuW})}
\begin{itemize}
  \item[(a)] %Let $T\sim\bP_{q\otimes\kappa}^{\rm col}$ be a coloured ${\rm GW}(q,\kappa)$-tree. Then
    For all Borel-measurable $k\colon (0,\infty)\rightarrow[0,\infty)$, $j\ge 2$, $\varepsilon_i\in\{0,1\}$ and
    $\cF^\bH$-measurable $f_i\colon \bT^\bH\rightarrow[0,\infty)$, $i=1,\ldots,j$, we have
     \begin{eqnarray*}\label{branching prop expec}
\lefteqn{\mathbb{E}_{q\otimes\kappa}^{p\rm-col}\left[k(\zeta_\emptyset)\prod_{i=1}^{j}f_i
(\overline{\theta}_i); \nu_\emptyset=j;(\gamma_1,\ldots,\gamma_j)=(\varepsilon_1,\ldots,\varepsilon_j)\right]}\nonumber\\
&
&=\int_{(0,\infty)}k(z)\kappa(dz)\ q(j)\ g(p)^{j_r}(1-g(p))^{j_b}\ \prod^{j}_{i=1}\mathbb{E}_{q\otimes\kappa}^{p\rm-col}[f_i|\gamma_{\emptyset}=\varepsilon_i]
     \end{eqnarray*}
    where $j_r=\varepsilon_1+\cdots+\varepsilon_j$ and $j_b=j-j_r$ are the numbers of red and black colour marks.
  \item[(b)] %Let $T\sim\bP_{q\otimes\kappa}^{\rm col}$ be a coloured ${\rm GW}(q,\kappa)$-tree and
    %$t\ge 0$. Then f
    For all $t\ge 0$ and $\cF^\bH$-measurable
    $f_u\colon \bT^\bH\rightarrow[0,\infty)$  $$\bE_{q\otimes\kappa}^{p\rm-col}\left[\left.\prod_{u\in\overline{\pi}_t}f_u(\overline{\theta}_{u,t})\right|\cF_t^\bH\right]=\prod_{u\in\overline{\pi}_t}\left.\bE_{q\otimes\kappa}^{p\rm-col}[f_u(\overline{\theta}_{\emptyset,s})|\zeta_\emptyset>s]\right|_{s=t-\alpha_u}$$
    In the exponential case $\kappa={\rm Exp}(c)$, this simplifies to
    \begin{equation}\bE_{q\otimes {\rm Exp}(c)}^{p\rm-col}\left[\left.\prod_{u\in\overline{\pi}_t}f_u(\overline{\theta}_{u,t})\right|\cF_t^\bH\right]=\prod_{u\in\overline{\pi}_t}\bE_{q\otimes{\rm Exp}(c)}^{p\rm-col}[f_u]
    \label{bpexp}
    \end{equation}
\end{itemize}
%on for each leaf of a random tree $\tau$, we mark it
%randomly and independently,  ``black or a mark $\varepsilon_u=1$" with
%probability $1-p$, and ``red or a mark $\varepsilon_u=0$" with probability
%$p$. If all the leaves are red, then we colour the whole tree,
%including the progenitor $\emptyset$ in red. A black vertex has a black
%``parent". If at least one leaf $u$ is black, then we will have a
%tree with black progenitor.
\paragraph{Reduction procedure to identify the ``black tree'' in a two-colours tree}
\begin{itemize}\item We can extract $\widetilde{T}^{p \rm -rdc}_{\rm sub}=\{(u,\zeta_u)\in\bU\times(0,\infty)\colon (u,\zeta_u,1)\in T^{p\rm-col}\}$, the individuals of $T^{p\rm-col}$ with black colour marks. If $\widetilde{T}^{p \rm -rdc}_{\rm sub}\neq\varnothing$,
we rename the individuals of $\widetilde{T}^{p\rm-rdc}_{\rm sub}$ by the unique injection
$$\iota\colon \widetilde{\tau}_{\rm sub}^{p \rm -rdc}=\{u\in\bU\colon (u,\zeta_u)\in\widetilde{T}_u^{p\rm-rdc}\}\rightarrow\bU,$$
that is increasing for the lexicographical total order on $\bU$, maps onto an element
$\tau^{p \rm -rdc}_{\rm sub}$ of $\bT$ and is compatible with the genealogical partial orders. We refer to the image tree $T^{p\rm-rdc}_{\rm sub}=\{(\iota(u),\zeta_{u})\colon u\in\widetilde{\tau}_{\rm sub}^{p\rm-rdc}\}$ as the \em $p$-reduced subtree of $T$\em.
\item As a further reduction, we remove single-child individuals and add their lifetimes to the child's lifetime. Formally, we define
  $\widetilde{\tau}^{p\rm-rdc}=\{v \in \tau^{p \rm -rdc}_{\rm sub}\colon {\nu_{v}}(\tau^{p \rm -rdc}_{\rm sub}) \neq 1\}$,
and $$\widetilde{\zeta}_u=\sum_{i=J_u}^n\zeta_{u_1\cdots u_i}(\tau^{p\rm-rdc}_{\rm sub}),
\mbox{ where }J_u=\sup\{j\colon \nu_{u_1\cdots u_i}(\tau^{p\rm-rdc}_{\rm sub})\!=\!1\mbox{ for all $i\in\{j,\ldots,n-1\}$}\},$$
for all $u=u_1\cdots u_n\in\widetilde{\tau}^{p\rm-rdc}$.
%this can be done representing $\tau_{\rm sub}^{p\rm-rdc}$ as a
%graph-theoretic tree $(V_{\rm sub},E_{\rm sub})$ with vertex set $V_{\rm sub}$ and edge set
%$E_{\rm sub}$ as follows:
%$$V=\{\emptyset\}\cup \{v \in \tau^{p \rm -rdc}_{\rm sub}:
%{\nu_{v}}(\tau^{p \rm -rdc}_{\rm sub}) \neq 1\}$$ and
%$$E=\{(u\rightarrow uj):uj\in\tau^{p\rm-rdc}_{\rm sub}\textrm{Merge of the edges between $u$ %and $v$ exluding  $u$ and $v$,} \quad u,v \in V.\}$$
%Put on $V$ the partial order inherited from $\tau^{p \rm -rdc}_{\rm
%sub}$.
Again, there is a unique injection $\iota^\prime\colon \widetilde{\tau}^{p\rm-rdc}\rightarrow\bU$ that is
increasing for the lexicographical total order on $\bU$, maps onto an element
$\tau^{p\rm-rdc}$ of $\bT$ and is compatible with the genealogical partial orders. We refer to the image tree $T^{p\rm-rdc}=\{(\iota^\prime(u),\widetilde{\zeta}_u)\colon u\in\widetilde{\tau}^{p\rm-rdc}\}$ as the \em $p$-reduced tree \em (or as the \em black tree\em).
\end{itemize}
%Then $(V,E)$ and the $\emptyset$ is an ordered tree in ``one
%to one'' correspondence to a unique element in $\mathbb{T}^{\rm
%%discr}$, denoted $\tau^{p \rm -rdc}$. It is not hard to show that the
%black subtree $\tau^{p \rm -rdc}_{\rm sub}$ and the black tree $\tau^{p
%\rm-rdc}$ are both Galton-Watson trees.
Figure \ref{discperc} in the Introduction illustrates the reduction procedure. 
\begin{rem}\label{remdet}\rm \begin{enumerate}\item[(a)] The reduction procedure is transitive in that for independent colouring, we have $(T^{(1-\overline{p}_1)\rm-rdc})^{(1-\overline{p}_2)\rm-rdc}\overset{{\rm (d)}}{=}  T^{(1-\overline{p}_1\overline{p}_2)\rm-rdc}$. In particular, colouring for $T^{(1-\overline{p}_1)\rm-rdc}$ and
$T^{(1-\overline{p}_3)\rm-rdc}$ for $\overline{p}_3<\overline{p}_1$ can be coupled such that $T^{(1-\overline{p}_3)\rm-rdc}=(T^{(1-\overline{p}_1)\rm-rdc})^{(1-\overline{p}_3/\overline{p}_1)\rm-rdc}$.
\item[(b)] Although we have used notation for a random ${\rm GW}(q,\kappa)$-tree $T^{p\rm-col}$ with leaves
coloured independently, note that the reduction of $T^{p\rm-col}$ to a black tree is a purely deterministic procedure. Our focus here has been on the technical framework and how it is used to formulate relevant examples. We postpone to Section \ref{Growth of GW exponential edge lengths} the review of further developments, notably of the reconstruction/growth procedures that reverse the reduction.
\end{enumerate}
\end{rem}
%------------------------------------------------------------------------------------------------------------------discrete trees with immigration
\subsubsection{Bushes and forests -- models with several progenitors and immigration}
\label{discr BP with immigration section} Branching processes
with immigration have been studied widely (see e.g.  Athreya and Ney 
\cite{AtN} and Jagers \cite{Jag-68}). We consider the model, where
immigrants arrive at the times $S_i$, $i\ge 1$, of a renewal
process $J_t=\#\{i\ge 1\colon S_i\le t\}$, i.e. where $S_0=0$ and the
interarrival times $S_i-S_{i-1}$, $i\ge 1$, are
independent and identically distributed random variables with a common
distribution on $(0,\infty)$ that we denote by $\chi$. At each
immigration time $S_i$ the number $N_i$ of immigrants is independent and
has a common distribution $\eta$ on $\bN$.
Each immigrant produces offspring independently according to the rules
of ${\rm GW}(q,\kappa)$-trees.

Denoting by $Z^{(i)}_{t-S_i}$ the size at time $t$ of the population of
immigrants arriving at time $S_i$,
\begin{equation}\label{discrimm} Y_t=\sum^{J_t}_{i=1}Z^{(i)}_{t-S_i}\end{equation}
is the total population size at time $t\ge 0$. Here, $Z^{(i)}$ are independent sums of $N_i$
independent Bellman-Harris processes with offspring and lifetime distributions $q$ and $\kappa$ as in Section \ref{Lifetime marks and discr BP in cont. time}.

To capture the genealogical trees of the population, we use the notion of a \em bush \em as a 
random sequence $B=(T_{(1)},\ldots,T_{(N)})$ of independent trees, and the notion of a \em
forest \em as a point process $F=(B(t),t\ge 0)$ of independent bushes
$$B(S_i)=B^{(i)},\quad i\ge 0,\qquad B(t)=\partial,\quad t\not\in\{S_i,i\ge 1\}\quad \mbox{for a cemetery point $\partial$.}$$
\paragraph{${\rm GW}(q,\kappa,\beta)$-bushes and ${\rm GWI}(q,\kappa,\eta,\chi)$-forests}
\begin{itemize}\item A ${\rm GW}(q,\kappa,\beta)$-bush is a bush $B=(T_{(1)},\ldots,T_{(N)})$ of
    independent ${\rm GW}(q,\kappa)$-trees $T_{(j)}$, where $N\sim{\rm Poi}(\beta)$.
  \item A ${\rm GWI}(q,\kappa,\eta,\chi)$-forest is a forest $F=(B(t),t\ge 0)$ where each
    bush $B(S_i)=B^{(i)}$ is associated with immigration at the times $S_i$ of a renewal
    process with inter-arrival distribution $\chi$ and consists of an independent
    $\eta$-distributed number $N_i$ of trees $T_{(j)}^{(i)}$.
\end{itemize}

It is straightforward to transfer the notions of colouring and reduction to the setting of bushes
and forests, since they apply tree by tree. We will slightly abuse notation and write $u\in B$ to refer
to individuals in a bush, instead of writing formally $u=(i,u^{\prime\prime})$ with $u^{\prime\prime}\in T_{(i)}$. Similarly, $u\in F$ means $u=(t,u^{\prime})$ with $u^\prime\in B(t)$ in the sense just defined. We will also abuse notation
$\nu_u$, $\zeta_u$ and $\gamma_u$ accordingly for $u\in B$ and $u\in F$.

%Moreover, if we, at time $0$, start such immigration process
%simultaneously with a GW process, which has the same offspring
%distribution $q$, the at time $t$ the size $Z_t$ of offspring alive
%is $$Z_t=Y_t + I_t.$$ We denote ${\rm GWI}(q,\kappa, \eta, \gamma)$ for a
%Galton-Watson process with offspring distribution $q$, i.i.d.
%edge-length with distribution $\kappa,$ i.i.d. immigrations with
%distribution $\eta$ and i.i.d. inter arrival times of immigration
%with distribution $\gamma.$

%-------------------------------------------------------------------------------------------------------------------CSBP and CBI intro
\subsection{Continuous-state branching processes and immigration}\label{contstate}
\label{CSBP and CBI} We have looked at branching processes with
immigration in the discrete state space $\bN_0=\{0,1,2,\ldots \}$ in
continuous time. In this section we will recall (Markovian) continuous-state
branching processes with immigration, where the state space
(population size) will no longer be $\mathbb{N}_0$ but
$[0,\infty)$, and time is also continuous.

%--------------------------------------------------------------------------------------------------------------------subordinators and g.i.d.
\subsubsection{Subordinators and geometric infinite divisibility}
\label{sub and g.i.d section}
%Subordinators are right-continuous increasing process with independent
%and homogeneous increments.
\begin{defi}\label{Subordinator}
\rm An increasing right-continuous stochastic process $\sigma=(\sigma(t), t \geq 0)$ in $[0,\infty)$ is called a \em subordinator \em if it has stationary independent increments, i.e. if for every
$u, t\ge 0$, the \em increment \em $\sigma(t+u)-\sigma(t)$ is independent of $(\sigma(s),s\le t)$ and $\sigma(t+u)-\sigma(t)\overset{{\rm (d)}}{=} \sigma(u)$.
\end{defi}
The distribution of a subordinator $\sigma$ on the space $\bD([0,\infty),[0,\infty))$ of functions $f\colon [0,\infty)\rightarrow[0,\infty)$ that are right-continuous and have left limits equipped with the Borel sigma-algebra generated by Skorohod's topology, see e.g. \cite{JaS}, is specified by the Laplace transforms of
its one-dimensional distributions. For every $t\ge 0$ and $r \ge 0$,
\begin{equation}\label{LaExponent}
\mathbb{E}(\exp\{-r \sigma(t)\})= \exp\{-t \phi (r)\}
\end{equation}
where the function $\phi \colon  [0, \infty) \rightarrow [0, \infty)$ is
called the \it Laplace exponent \rm of $\sigma$.
There exist a unique real number $d\ge 0$ and a
unique measure $\Lambda$ on $(0, \infty)$ with $\int (1 \wedge x) \Lambda
(dx) < \infty$, such that for every $r \geq 0$
\begin{equation}\label{LK1}
\phi(r)= dr + \int_{(0,\infty)}(1-e^{-r x}) \Lambda
(dx).
\end{equation}
Conversely, any function $\phi$ of the form (\ref{LK1}) is the
Laplace exponent of a subordinator, which can be constructed as $(dt+\sum_{s\le t}\Delta\sigma_s,t\ge 0)$ for a Poisson point
process $(\Delta\sigma_s,s\ge 0)$ in $(0,\infty)$ with intensity measure $\Lambda$. Equation (\ref{LK1}) is
referred to as the L\'evy-Khintchine representation of $\phi$.
We refer to Bertoin
\cite{Ber-Sub} for an introduction to subordinators and their
applications.

\begin{defi}
\rm A random variable $X$ is
\em geometrically infinitely divisible \em (g.i.d.) if for all
$\alpha\in(0,1)$
\begin{equation}\label{GEO}
X \overset{{\rm (d)}}{=}  \sum_{j=1}^{G^{(\alpha)}}X_{\alpha}^{(j)}, \qquad j\ge 1,
\end{equation}
for a sequence $X_\alpha^{(j)}$, $j\ge 1$, of independent identically distributed (i.i.d.) random variables and an independent $G^{(\alpha)}\sim{\rm Geo}(\alpha)$:
$$
\mathbb{P}(G^{(\alpha)}=k)=\alpha(1-\alpha)^{k-1},\qquad
k\ge 1.
$$
\end{defi}
For example, an exponential random variable $X \sim {\rm Exp}(c)$
is g.i.d.\ since (\ref{GEO}) holds for $X_{\alpha}^{(j)} \sim {\rm Exp}(c/\alpha )$. The class of g.i.d.\ distributions can be characterised as follows.
\begin{lemm}[\cite{KMM}] \label{sub and geo}
A random variable $X$ is g.i.d.\ if and
only if it can be expressed as $X \overset{{\rm (d)}}{=}  \sigma(V)$, where $\sigma=(\sigma(t), t \geq 0)$ is
a subordinator and $V \sim {\rm Exp}
(c)$ independent of $\sigma$ for one equivalently all $c\in(0,\infty)$.
\end{lemm}
\noindent Indeed we can then express ($V$ and hence) $X$ as
$$V=\sum_{j=1}^{G^{(\alpha)}}V_{\alpha}^{(j)}\qquad\mbox{and}\qquad X=\sigma_V=\sum_{j=1}^{G^{(\alpha)}}\left(\sigma\left(V_\alpha^{(1)}+\cdots+V_\alpha^{(j)}\right)-\sigma\left(V_\alpha^{(1)}+\cdots+V_\alpha^{(j-1)}\right)\right).$$
\begin{lemm}\label{uniqgeodiv} If $X$ is such that {\rm(\ref{GEO})} holds for some $\alpha\in(0,1)$, the distribution of $X_\alpha^{(j)}$ is unique.
\end{lemm}
\begin{pf}  $\displaystyle\bE(e^{-rX})=\frac{\bE(e^{-rX_\alpha^{(1)}})\alpha}{1-(1-\alpha)\bE(e^{-rX_\alpha^{(1)}})}\quad\Rightarrow\quad\bE(e^{-rX_{\alpha}^{(1)}})=\frac{\bE(e^{-rX})}{\alpha+(1-\alpha)\bE(e^{-rX})}.$
\end{pf}
%-------------------------------------------------------------------------------------------------------------------------------CSBP and CBI explanation
\subsubsection{Continuous-state branching processes}\label{CSBP and CBI explanation
section} Continuous-state branching processes were
introduced by Jirina \cite{Jir} and Lamperti \cite{Lam-CSBP}.
They are the limiting processes of sequences of rescaled
Galton-Watson processes as the initial population size tends
to infinity and the mean lifetime tends to zero. In this section we follow Le Gall
\cite{LeG-snake}.
\begin{defi}\rm
A continuous-state branching process (CSBP) is a right-continuous Markov
process $(Z_t, t \geq 0)$ in $[0,\infty)$, whose transition
kernels $P_t(x, dz)$ are such that for every $t\ge 0$, $x \geq 0$ and $x^\prime\ge 0$ we have
$$P_t(x+x^\prime,\cdot)=P_t(x,\cdot)*P_t(x^\prime,\cdot),$$
where $*$ denotes convolution. In other words, if for a given transition kernel we denote, for each $x\ge 0$, by $Z^x$ a CSBP starting from $Z_0^x=x$, then for $\widetilde{Z}^{x^\prime}\overset{{\rm (d)}}{=}Z^{x^\prime}$ independent, we require $Z^x_t +\widetilde{Z}^{x^\prime}_t\overset{{\rm (d)}}{=}  Z^{x+x^\prime}_t$.
\end{defi}

The transition kernel is specified by the Laplace transforms
\begin{equation}\label{CSBP and psi solution u}
    \mathbb{E}(\exp\{-r Z_t\}|Z_0=x) =\exp\{-x u_t(r)\}
\end{equation}
where $u_t\colon [0,\infty)\rightarrow[0,\infty)$. In fact, $(t,r)\mapsto u_t(r)$
is necessarily the unique non-negative solution of
\begin{equation}\label{u solution of psi}
    u_t(r)+\int^t_0 \psi (u_s (r)) ds = r \qquad
    \mbox{or} \qquad \frac{\partial u_t(r)}{\partial t} = -
    \psi (u_t (r))
\end{equation}
with $u_0 (r)=r$, for some $\psi\colon [0,\infty)\rightarrow\bR$ of the form 
\begin{equation}\label{LK2}
\psi(r)=br+a r^{2} + \int_{(0,
\infty)}(e^{-r x}-1+ r x \mathbf{1}_{\{x<1\}})\Pi(dx),%\quad b \in \mathbb{R},\ a \geq 0
\end{equation}
where $b \in \mathbb{R}$, $a \geq 0$ and $\Pi$ is a measure on $(0,\infty)$ with $\int(1 \wedge x^2)\Pi(dx)<\infty$ and where $\psi$ satisfies the non-explosion condition $\int_{0+}|\psi(r)|^{-1}dr=\infty$, see \cite{Grey}. 
Equation (\ref{LK2}) is referred to as the L\'evy-Khintchine representation of $\psi$. The process $(Z_t,t\ge 0)$ is then called a CSBP \em with branching mechanism $\psi$\em, or a
${\rm CSBP}(\psi)$. We denote by $\bP^x_\psi$ the distribution of $(Z_t^x,t\ge 0)$ on $\bD([0,\infty),[0,\infty))$.  There also exists a sigma-finite measure $\Theta_\psi$ on $\bD([0,\infty),[0,\infty))$, such 
that 
$$(Z_t^x,t\ge 0)\overset{{\rm (d)}}{=} \left(\sum_{0\le y\le x}E_t(y),t\ge 0\right),\qquad x\ge 0,$$
where $(E(y),y\ge 0)$ is a Poisson point process in $\bD([0,\infty),[0,\infty))$ with intensity measure $\Theta_\psi$. 

\subsubsection{Continuous-state branching processes with immigration}

Similarly, a discrete-state branching process with immigration has a
continuous analogue, the continuous-state branching process with
immigration, CBI for short. Following Kawazu and Watanabe \cite{KaW-71}, see also 
\cite{Duq-imm,Lam-CBI}, besides the branching mechanism
$\psi$ for a CSBP, we also have an immigration mechanism $\phi$ of the form (\ref{LK1}) for
the CBI, which we then refer to as ${\rm CBI}(\psi,\phi)$.

A ${\rm CBI}(\psi,\phi)$ is a Markov process $(Y_t, t \geq 0)$ on $[0,\infty)$ whose transition kernels are
characterized by their Laplace transform, which in terms of $\phi$ and $u_t(r)$ as in (\ref{u solution of psi}) satisfy
%\begin{equation}\label{CBI transition LT form}
$$    \mathbb{E}(\exp \{-r Y_t \}|Y_0=x)=\exp\left\{-x u_t(r)-\int^t_0 \phi(u_s(r))ds\right\}.
$$%\end{equation}
In fact, a subordinator
$\sigma=(\sigma(t),t\ge 0)$ with Laplace exponent $\phi$ can be viewed as a pure immigration process ${\rm CBI}(\phi,0)$. Indeed, a general ${\rm CBI}(\psi,\phi)$ is such that by time $t\ge 0$ a population of total size $\sigma(t)$ has immigrated and evolved like a ${\rm CSBP}(\psi)$; specifically, consider a Poisson
point process $(E^s,s\ge 0)$ in $\bD([0,\infty),[0,\infty))$ with intensity measure $d\Theta_\psi+\int_{(0,\infty)}\bP^x_\psi\Lambda(dx)$, 
then in analogy with (\ref{discrimm})
$$ Y_t=\sum_{s\le t}E_{t-s}^s,\qquad t\ge 0,
$$
is a ${\rm CBI}(\psi,\phi)$. This follows by the exponential formula and properties of Poisson point processes.

Examples of continuous-state branching processes with immigration include (sub)critical CSBP conditioned on survival, see \cite{Lam-07} and literature therein.

%----------------------------------------------------------------------------------------------------------------------Growth of GW trees bush Exponential edge lengths
\subsection{Growth of Galton-Watson bushes with exponential edge
lengths}\label{Growth of GW exponential edge lengths} In Theorem \ref{thm1}, we consider
families of ${\rm GW}(q_\lambda,{\rm Exp}(c_\lambda),\beta_\lambda)$-bushes $B_\lambda$,
$\lambda\ge 0$, with the
consistency property that any two bushes, for parameters $\mu<\lambda$ say, are related by $p$-reduction as formally defined in Section \ref{coloured leaves and 2-colour BP}, for $p=1-\mu/\lambda$. The choice of $p$ is dictated (up to a positive power for the ratio) by the consistency requirement that the relation holds for \em all \em $\mu$ and \em all \em
$\lambda$ (cf.\ Remark \ref{remdet}(a)).

The equivalence of (i) and (ii) in Theorem \ref{thm1} is a statement purely about offspring distributions $(q_\lambda,\lambda\ge 0)$. The reason for including the other two parameter sequences $(c_\lambda,\lambda\ge 0)$ and $(\beta_\lambda,\lambda\ge 0)$ in the remainder of the statement is simplicity. Specifically, if we consider trees without lifetimes and
hence without embedding in time, the removal of single-child individuals will be more artificial as it reduces the height of the trees; if we look at trees instead of bushes, the reduced tree will only be GW if we condition on the existence of at least one black leaf, and this does not lead to consistent families of random trees on the same probability space.

In \cite{DuW}, the bushes of Theorem \ref{thm1} are used to construct
L\'evy trees as a strong representation of the genealogy of the limiting ${\rm CSBP}(\psi)$ under some extra conditions on $\psi$. In a weaker sense, the family $(B_\lambda,\lambda\ge 0)$
itself is already a representation of the genealogy of ${\rm CSBP}(\psi)$, under
no conditions on $\psi$ other than $\psi(\infty)=\infty$ to exclude the case of increasing CSBPs corresponding to ``no death'' i.e.\ ``no leaves''. Roughly, $B_\lambda$ is the genealogy of a Poisson sample chosen among
all individuals with intensity proportional to $\lambda$; as $\lambda\rightarrow\infty$, the set
of individuals included becomes dense.

%-------------------------------------------------------------------------------------------------------------Consistent families of GW(q,c,beta)-bushes
\subsubsection{Consistent families of ${\rm GW}(q_\lambda, {\rm Exp}(c_\lambda), \beta_\lambda
)$-bushes}\label{Consistent family}

In \cite{DuW}, the parameters in a consistent family of
Galton-Watson bushes are represented in terms of a branching
mechanism $\psi$ that is just required to satisfy $\psi(\infty)=\infty$, so that $\psi$ is eventually increasing and has a right inverse $\psi^{-1}\colon[0,\infty)\rightarrow[\psi^{-1}(0),\infty)$:
\begin{equation}\label{offspr}\varphi_{q_\lambda}(s)=s+\frac{\psi(\psi^{-1}(\lambda)(1-s))}{\psi^{-1}(\lambda)\psi^\prime(\psi^{-1}(\lambda))},\qquad c_\lambda=\psi^\prime(\psi^{-1}(\lambda)),\qquad\beta_\lambda=\beta\psi^{-1}(\lambda).
\end{equation}
It follows from the derivation there that this $\psi$ and $(\widetilde{\psi},c)$ in the statement of Theorem \ref{thm1} are related in a linear way as
\begin{equation}\label{psipsitilde}\psi(r)=k_1\widetilde{\psi}(k_2 r),
\end{equation}
where $k_1=1/\widetilde{\psi}(1)=\psi^{-1}(1)\psi^\prime(\psi^{-1}(1))$ and $k_2=c\widetilde{\psi}(1)=1/\psi^{-1}(1)$. For the underlying CSBPs, the relationship (\ref{psipsitilde}) just means $Z_t=k_2\widetilde{Z}_{k_1k_2t}$, so $\psi$ and $\widetilde{\psi}$
essentially refer to the same CSBP. With this parameterisation, (\ref{gp}) can be expressed
more explicitly for $q=q_\lambda$ and $p=1-\mu/\lambda$ as
\begin{equation}\label{g in psi}
    g_\lambda(1-\mu/\lambda)=\mathbb{P}_{q_\lambda\otimes{\rm Exp}(c_\lambda)}^{(1-\mu/\lambda)\rm-col}(\gamma_{\emptyset}=1)
=1-\frac{\widetilde{\psi}^{-1}(\mu)}{\widetilde{\psi}^{-1}(\lambda)}=1-\frac{\psi^{-1}(\mu)}{\psi^{-1}(\lambda)}.
\end{equation}

%----------------------------------------------------------------------------------------------------------Reconstruction procedure for GW(q, c, beta)-bushes
\subsubsection{Analysis of the reduction procedure for ${\rm GW}(q_\lambda,{\rm Exp}(c_\lambda),
\beta_\lambda)$-bushes}\label{Recon for GW expo bushes1} In this section
we will study some key points in the reduction procedure in the setting of Theorem \ref{thm1}.
These will be important for the proofs of Theorems \ref{thm2} and \ref{thm3}.
There are three steps in the reduction procedure from $\lambda$ to $\mu<\lambda$: colouring
with $p=1-\mu/\lambda$, passage to the $p$-reduced sub-bush and passage to the $p$-reduced bush.

In the last step, the lifetime $\zeta^{(1-\mu/\lambda) \rm -rdc}_{u}$ of individual $u\in T^{(1-\mu/\lambda) \rm -rdc}_\lambda$ is obtained combining a
number $G^{(\alpha)}_u$ of lifetimes of $T_\lambda$ when removing the
single-child individuals. By the two-colours branching property in Section \ref{coloured leaves and 2-colour BP}, the random numbers $G^{(\alpha)}_u$ are
${\rm Geo}(\alpha)$, where $\alpha$ is the probability that an individual in
$T^{(1-\mu/\lambda) \rm -rdc}_{\lambda, \rm sub}$ produce zero or more than two children. We can express $\alpha$ in terms of $\psi$:
\begin{lemm} \label{alpha and psi} Let $0\le\mu<\lambda<\infty$. Given that $u\in B^{(1-\mu/\lambda)\rm-rdc}_\lambda$, we have
  $G^{(\alpha)}_u\sim{\rm Geo}(\alpha)$ with
\begin{equation}\label{imporesult}
\alpha= \mathbb{P}\left(\left.\nu_{\emptyset}\left(T_{\lambda,\rm sub}^{(1-\mu/\lambda)\rm-rdc}\right)\neq 1\;\right|\,
\gamma_{\emptyset}\left(T_\lambda^{(1-\mu/\lambda)\rm-col}\right)=0\right)=\frac{\psi'(\psi^{-1}(\mu))}{\psi'(\psi^{-1}(\lambda))}.
\end{equation}
%where $\tau_\lambda^{(1-\mu/\lambda) \rm -rdc}(.)$ and
%$\tau_\lambda(.)$ are the generating functions of offspring
%distribution $q_\lambda^{(1-\mu/\lambda) \rm -rdc}$ and $q_\lambda$
%respectively.
\end{lemm}
\begin{pf} According to the definition of $\alpha$ and the two-colours branching property of Section \ref{coloured leaves and 2-colour BP},
\begin{eqnarray*}
1-\alpha&=&\mathbb{P}\left(\left.\nu_\emptyset\left(T_{\lambda,\rm sub}^{(1-\mu/\lambda)\rm-rdc}\right)=1\; \right|\, \gamma_\emptyset\left(T_\lambda^{(1-\mu/\lambda)\rm-col}\right)=0\right)\\
  &=&\frac{1}{1-g_\lambda(1-\mu/\lambda)}\sum^{\infty}_{j=2}{j\choose 1}q_\lambda(j)g_\lambda(1-\mu/\lambda)^{j-1}(1-g_\lambda(1-\mu/\lambda))\\
                                                    &=&\varphi_{q_\lambda}^\prime(g_\lambda(1-\mu/\lambda))
\end{eqnarray*}
and by (\ref{offspr}) and (\ref{g in psi}) we obtain
%\begin{equation}\label{alpha mu lambda}
$1-\alpha=1-\psi^\prime(\psi^{-1}(\mu))/\psi^\prime(\psi^{-1}(\lambda))$.
%\end{equation}
\end{pf}
In the reconstruction procedure reversing the reduction procedure, we will therefore subdivide each lifetime in the
${\rm GW}(q_\mu, {\rm Exp}(c_\mu),\beta_\mu)$-bush into a geometric random
number of $G^{(\alpha)}_u$ parts.

\begin{rem}\label{alphakappa}\rm By the branching property of coloured ${\rm GW}(q,\kappa)$-trees in Section \ref{coloured leaves and 2-colour BP}, lifetime marks are
  independent of colour marks. Therefore, Lemma \ref{alpha and psi} also holds for general $B_\lambda\sim{\rm GW}(q_\lambda,\kappa_\lambda,\beta_\lambda)$-bushes, not just for ${\rm GW}(q_\lambda,{\rm Exp}(c_\lambda),\beta_\lambda)$-bushes.
\end{rem}

In the case of ${\rm Exp}(c_\mu)$ lifetimes, given the lifetime $\zeta_u(T_\mu)$ the
conditional distribution of the random variable $N_u=G^{(\alpha)}_u-1$ follows a Poisson distribution:
\begin{prop}\label{Poi&Exp}
Let $\zeta$ be a random variable having distribution $\zeta \sim
{\rm Exp}(\alpha c)$. Suppose that $\zeta$ is subdivided into $G^{(\alpha)}$
independent parts, $\zeta=\zeta_1+\cdots+\zeta_{G^{(\alpha)}}$, where
$\zeta_i \sim {\rm Exp}(c)$ and $G^{(\alpha)} \sim {\rm Geo}(\alpha)$ are independent. Then given $\zeta=z$ for $z \geq 0$, we have
\begin{equation}\label{Poi with exp condition}
\mathbb{P}(G^{(\alpha)}=k \ | \ \zeta=z )
=\frac{((1-\alpha)cz)^{k-1}e^{- (1-\alpha) c z}}{(k-1)!}.
\end{equation}
\end{prop}
\begin{pf} This is, of course, well-known in the context of Poisson processes, but let us
give a direct argument and write the left hand side as a conditional
expectation
$\mathbb{P}(G^{(\alpha)}=k \ | \ \zeta=z )=\mathbb{E}(\mathbf{1}_{\{G^{(\alpha)}=k\}} |
\zeta=z)$. We also set $g_k(z)={((1-\alpha)cz)^{k-1}e^{- (1-\alpha) c
z}}/{(k-1)!}$.

$$\mbox{\bf Claim: \rm $\mathbb{E}(f(\zeta)g_k(\zeta))=\mathbb{E}(f(\zeta)\mathbf{1}_{\{G^{(\alpha)}=k\}})$ for all measurable $f \geq 0$.}$$
As $\zeta \sim {\rm Exp}(\alpha c)$,
$\displaystyle%\begin{eqnarray*}
\mathbb{E}(f(\zeta)g_k(\zeta)) = \int_{0}^{\infty}f(z)g_k(z)\alpha c
e^{-
\alpha c z}dz
%&=& \int_{0}^{\infty} f(z)\frac{((1-\alpha)c z)^{n-1}e^{- (1-\alpha)
%c z}}{(n-1)!} \alpha c e^{- \alpha c z}dz\\
= \int_{0}^{\infty}f(z)\frac{\alpha c ((1-\alpha) c
z)^{k-1}e^{-cz}}{(k-1)!}dz
$. %\end{eqnarray*}
On the other hand, since $\zeta=\zeta_1+\cdots+\zeta_{G^{(\alpha)}}$, the conditional
distribution of $\zeta$ given $G^{(\alpha)}=k$ is Gamma$(c, k)$. Therefore,
$$%\begin{eqnarray*}
\mathbb{E}(f(\zeta)\mathbf{1}_{\{G^{(\alpha)}=k\}}) =
\mathbb{P}(G^{(\alpha)}=k)\int_{0}^{\infty}f(z)\frac{z^{k-1}c^k
e^{-cz}}{(k-1)!}dz
= \int_{0}^{\infty}f(z)\frac{\alpha c ((1-\alpha) c
z)^{k-1}e^{-cz}}{(k-1)!}dz,
$$%\end{eqnarray*}
and so $g_k(z)$ is a version of the conditional
probability $\mathbb{P}(G^{(\alpha)}=k \ | \ \zeta=z )$, as claimed.
\end{pf}%\begin{flushright}$\Box$\end{flushright}
%\begin{rem}\rm:
%\begin{enumerate}
%  \item \rm In Proposition \ref{Poi&Exp}, given the $\zeta=z \ (x>0),$  $N=Y^{(\alpha)}-1$ of has
%a Poisson distribution with parameter $(1-\alpha)c l.$
%  \item In terms of a black tree with an edge of length $\zeta=z$, and $\zeta \sim {\rm Exp} (c_\mu=\overline{\psi}'(\overline{\psi}^{-1}(\mu)) ),
%  $ $N$ is the number of
%vertices added to subdivide $\zeta$ into $Y^{(\alpha)}=N+1$ parts in
%the reconstruction procedure. Applying Proposition \ref{Poi&Exp},
%from the Poisson process on $[0,l],$ we know where to subdivide this
% $\zeta$-long edge with $N$ vertices.
%\end{enumerate}
%\end{rem}
Moreover, we can also find the conditional joint distribution of the
lifetimes $(\zeta_1,\ldots,\zeta_{G^{(\alpha)}-1})$ given $\zeta=z$ and $G^{(\alpha)}=k$ as follows:
\begin{equation}\label{joint exponential}
f_{\zeta_1,\zeta_2,\ldots,\zeta_{k-1}|G^{(\alpha)}=k,\zeta=z} (y_1,
y_2,\ldots,y_{k-1})%=\frac{\psi^\prime(\psi^{-1}(\mu))^k
%e^{-\psi^\prime(\psi^{-1}(\mu))l}(k-1)!}{\psi^\prime(\psi^{-1}(\mu))^k
%l^{k-1}e^{-\psi^\prime(\psi^{-1}(\mu))l}}
=\frac{(k-1)!}{z^{k-1}}
\end{equation}
for all $y_1>0,\ldots,y_{k-1}>0$ such that $y_1+\cdots+y_{k-1}<z$.
%\newpage

In the middle step of the reduction procedure, all red individuals are removed. \cite{DuW} noted that it is a
consequence of the two-colours branching property, see Section \ref{coloured leaves and 2-colour BP} here, that they form independent ${\rm GW}(q_{\lambda,\rm red}^{(1-\mu/\lambda)\rm-col},{\rm Exp}(c_\lambda))$-trees, ``red trees'', where
%--------------%%%%%%%%-----------%%%%%%%%----------%%%%%%------offspring distributions-%%%%%--------------------%%%%%%%%%%%---------------------------
\begin{equation}\label{red offspring distribution}
\varphi_{q_{\lambda,\rm red}^{(1-\mu/\lambda)\rm-col}}(s)=s+\frac{\psi_\mu(\psi_\mu^{-1}(\lambda-\mu)(1-s))}{\psi_\mu(\lambda-\mu)\psi_\mu^\prime(\psi_\mu^{-1}(\lambda-\mu))},\qquad\mbox{with $\psi_\mu(r)=\psi(\psi^{-1}(\mu)+r)-\mu$.}
%(m)=\left\{\begin{array}{ll}
%                     \frac{|\psi^{(m)}(\psi^{-1}(\lambda))|}{m!\psi'(\psi^{-1}(\lambda))(\lambda)-\psi^{-1}(\mu))^{m-1)}} & \mbox{if $m \geq 1$;}\\
%                     \frac{\psi(\psi^{-1}(\lambda))(\lambda-\mu)}{\lambda (\psi^{-1}(\lambda)-\psi^{-1}(\mu))(\psi'(\psi^{-1}(\lambda)))}     & \mbox{if  $m=0$.}
%               \end{array}
%        \right.
\end{equation}
Furthermore, the numbers of red trees removed at the branchpoints are conditionally independent, and
given that a branchpoint has $m \geq 1$ subtrees containing black leaves, the generating function of the number of red trees can be expressed in terms of the $m$th derivative $\psi_\mu^{(m)}$ of $\psi_\mu$:
\begin{equation}
\label{redgraft} \frac{\psi_\mu^{(m)}(\psi_\mu^{-1}(\lambda-\mu)(1-s))}{\psi_\mu^{(m)}(0)}\mbox{ for $m\ge 2$ or }\frac{\psi_\mu^\prime(\psi_\mu^{-1}(\lambda-\mu)(1-s))-\psi_\mu^\prime(\psi_\mu^{-1}(\lambda-\mu))}{\psi_\mu^\prime(\psi_\mu^{-1}(\lambda-\mu)-\psi_\mu^\prime(0)}\mbox{ for $m=1$.}
\end{equation}

%%%%%%%%%%%%%%%%%%%%%%%%%------------------------offspring distribution above--------------------------------%%%%%%%%%%%%%%%%%%%----------------------

\subsubsection{Reconstruction procedure for ${\rm GW}(q_\lambda,{\rm Exp}(c_\lambda),
\beta_\lambda)$-bushes}\label{Recon for GW expo bushes2}

Let $B_\mu=(T^{(1)}_\mu, T^{(2)}_\mu,\ldots ,T^{(N_\mu)}_\mu)$ be a ${\rm GW}(q_\mu,{\rm Exp}(c_\mu),\beta_\mu)$-bush. We will construct $B_\lambda$.
\begin{enumerate}\item[1.] For each individual $u$ in $T_\mu^{(i)}$, given the lifetime
  $\zeta_u^{(i)}=z$, subdivide into a random number $G_u^{(\alpha,i)}$ of parts with
  distribution (\ref{Poi with exp condition}), where $\alpha$ is as in (\ref{imporesult})
  and the parts $(\zeta_{u,1}^{(i)},\ldots,\zeta_{u,G_u^{(\alpha,i)}}^{(i)})$ have joint
  distribution (\ref{joint exponential}). Now for each $i\in\{1,\ldots,N_\mu\}$, there is a
  unique injection
  $$(\iota^\prime_i)^{-1}\colon \tau^{(i)}_\mu=\{u\in\bU\colon (u,\zeta_u^{(i)})\in T^{(i)}_\mu\}\rightarrow\bU$$
  such that $(\iota^\prime_i)^{-1}(\emptyset)=1^k$ with $k=G_\emptyset^{(\alpha,i)}-1$ and
  $(\iota^\prime_i)^{-1}(uj)=(\iota_i^\prime)^{-1}(u)j1^{k}$ with $k=G_{uj}^{(\alpha,i)}-1$, for all
  $uj\in\tau_\mu^{(i)}$, where $1^k$ is a string of $k$ letters $1$. We define
  $$\widehat{T}_{\mu}^{(i)}=\{(1^{n-1},\zeta_{\emptyset,n}^{(i)})\colon 1\!\le\! n\!\le\! G_\emptyset^{(\alpha,i)}\}\cup\{((\iota_i^\prime)^{-1}(u)j1^{n-1},\zeta_{u,n}^{(i)})\colon uj\in\tau_\mu^{(i)},1\!\le\! n\!\le\! G_u^{(\alpha,i)}\}.$$
\item[2.] For each individual $u$ in $\widehat{T}_{\mu}^{(i)}$, given
  $\widehat{\nu}_u^{(i)}=m \geq 1$ children,
  consider a random number $C_u^{(i)}$ of further children with distribution (\ref{redgraft}) and a
  uniform random permutation $\varrho_u^{(i)}$ among the $(m+C_u^{(i)})!/m!$ permutations with
  $\varrho_u^{(i)}(1)<\cdots<\varrho_u^{(i)}(m)$. Let $T^{(u,1,i)},\ldots,T^{(u,k,i)}$ with
  $k=C_u^{(i)}$ be independent
  ${\rm GW}(q_{\lambda,\rm red}^{(1-\mu/\lambda)\rm-rdc},{\rm Exp}(c_\lambda))$-trees with offspring
  distribution (\ref{red offspring distribution}). Then for each $i\in\{1,\ldots,N_\mu\}$, there
  is a unique injection
  $$(\iota_i)^{-1}\colon \widehat{\tau}^{(i)}_{\mu}=\{u\in\bU\colon (u,\widehat{\zeta}_u^{(i)})\in \widehat{T}_\mu^{(i)}\}\rightarrow\bU$$
  such that $(\iota_i)^{-1}(\emptyset)=\emptyset$ and $(\iota_i)^{-1}(uj)=(\iota_i)^{-1}(u)\varrho_u^{(i)}(j)$
  for all $uj\in\widehat{\tau}^{(i)}_\mu$. We define
    \begin{eqnarray*}\widehat{T}_\lambda^{(i)}&=&\{((\iota_i)^{-1}(u),\widehat{\zeta}_u^{(i)})\colon u\in\widehat{\tau}_\mu^{(i)}\}\\ &&\cup\{(\iota_i)^{-1}(u)\varrho_u(\widehat{\nu}_u^{(i)}+j)v,\zeta_v^{(u,j,i)})\colon u\in\widehat{\tau}_\mu^{(i)},1\le j\le C_u^{(i)},v\in\tau^{(u,j,i)}\}.
    \end{eqnarray*}
\item[3.] Given $N_\mu=n$, consider a random number
  $N_{\lambda}^{\rm red}\sim{\rm Poi}(\beta_\lambda-\beta_\mu)$ of further progenitors and a
  uniform random permutation $\varrho$ among the $(n+N_{\lambda}^{\rm red})!/n!$ permutations with
  $\varrho(1)<\cdots<\varrho(n)$. Let $\widehat{T}^{(n+1)}_\lambda,\ldots,\widehat{T}^{(n+k)}$
  with $k=N_{\lambda}^{\rm red}$ be independent
  ${\rm GW}(q_{\lambda,\rm red}^{(1-\mu/\lambda)\rm-rdc},{\rm Exp}(c_\lambda))$-trees with
  offspring distribution (\ref{red offspring distribution}). Then we finally define
  $$B_\lambda=(\widehat{T}_\lambda^{\varrho^{-1}(1)},\ldots,\widehat{T}_\lambda^{\varrho^{-1}(n+k)})\qquad\mbox{with $n=N_\mu$ and $k=N_\lambda^{\rm red}$.}$$
\end{enumerate}

\begin{rem}\label{Markov}\rm \begin{enumerate}\item[(a)] The constructed bush is indeed a ${\rm GW}(q_\lambda,{\rm Exp}(c_\lambda),\beta_\lambda)$-bush
and the pair obtained $(B_\mu,B_\lambda)$ has the same distribution as
$(B_\lambda^{(1-\mu/\lambda) \rm-rdc},B_\lambda)$. The intermediate trees
$\widehat{T}_\mu^{(i)}$ have the same distribution as the $(1-\mu/\lambda)$-reduced subtrees,
also jointly with the pair. If the split into two parts in the definition of $\widehat{T}_\lambda^{(i)}$ is used to assign black colour marks to the first part and red colour marks to the second part, then the resulting trees have the same distribution as the $(1-\mu/\lambda)$-coloured trees $B^{(1-\mu/\lambda)\rm-col}_\lambda$. We refer to Figure \ref{discperc} as a graphical illustration of the reduction and hence the reconstruction procedure.
%; Figure \ref{discperc} gives an example without lifetimes.
\item[(b)] In \cite{DuW}, this reconstruction procedure is also formulated for representations of the trees
in a space of $\bR$-trees, tree-like metric spaces that we briefly address in Section \ref{realtree}.
\item[(c)] It is a simple consequence of the reduction procedure and/or the reconstruction
procedure that $(B_\lambda,\lambda\ge 0)$ is an inhomogeneous \em Markov \em process in a suitable space of finite sequences of $(0,\infty)$-marked trees. Similarly, (\ref{bpexp}) can be strengthened in the present context to 
  $$\bE\left(\left.\prod_{u\in\overline{\pi}_t(B_\lambda)}f_u\left(\overline{\theta}_{u,t}\left(B_\lambda^{(1-\mu/\lambda)\rm-col}\right)\right)\right|\cG_{\lambda,t}\right)=\prod_{u\in\overline{\pi}_t(B_\lambda)}\bE\left(f_u\left(B_\lambda^{(1-\mu/\lambda)\rm-col}\right)\right),$$
  where $\cG_{\lambda,t}=\sigma\left\{\overline{\pi}_s(B_{\lambda^\prime}),\lambda^\prime \geq \lambda, s \leq t \right\}$, $\lambda\ge 0$, $t\ge 0$, and $B^{(1-\mu/\lambda)\rm-col}_\lambda$ is as in (a).
\end{enumerate}
\end{rem}

%-----------------------------------------------------------------------------------------------------------------limiting behaviour for \lambda recon
\subsubsection{Limiting behaviour as $\lambda \rightarrow \infty$} \label{limiting of lambda for edge} In the context of limiting results as $\lambda\rightarrow\infty$, we record the following corollary of Lemma \ref{alpha
and psi}.
\begin{coro}\label{9} If we fix $\mu>0$ in the setting of Lemma \ref{alpha and psi}, we have
$$\alpha=\alpha(\mu,\lambda)=\frac{\psi^\prime(\psi^{-1}(\mu))}{\psi^\prime(\psi^{-1}(\lambda))}\rightarrow\frac{\psi^\prime(\psi^{-1}(\mu))}{\psi^\prime(\infty)}\qquad\mbox{as $\lambda\rightarrow\infty$,}$$
where $\psi^\prime(\infty)$ means $\lim_{\lambda\rightarrow\infty}\psi^\prime(\lambda)$, in the following sense:
\begin{itemize}
  \item If $\psi^\prime(\infty)< \infty,$ then $\alpha(\mu,\lambda)\rightarrow\alpha_0(\mu)=\psi^\prime(\psi^{-1}(\mu))/\psi^\prime(\infty)$.
  \item If $\psi^\prime(\infty)=\infty$, then $\alpha(\mu,\lambda)\rightarrow\alpha_0(\mu)=0$.
\end{itemize}
\end{coro}

Note that this means that as $\lambda\rightarrow\infty$, lifetimes are cut into finite
${\rm Geo}(\alpha_0(\mu))$-distributed numbers of pieces in the first case, but into infinitely many
pieces in the second case.

%---------------------------------------------------------------------------------------------------------------------------------------------------
%---
%---------------------------------------------------------------------------------------------------------------------------------------------------
%---
%---------------------------------------------------------------------------------------------------------------------------------------------------
%---
%--------------------------------------From John----------------------------From John----------------------------------------------------------------
%---
%---------------------------------------------------------------------------------------------------------------------------------------------------
%---
%---------------------------------------------------------------------------------------------------------------------------------------------------
%---
%---------------------------------------------------------------------------------------------------------------------------------------------------
\begin{lemm} \label{CSBP Limit Lemma x}
Let $(Z_t^{\lambda},t\ge 0)$, $\lambda\ge 0$, be continuous-time Galton-Watson processes associated with a consistent family $(B_\lambda)_{\lambda\ge 0}$, of ${\rm GW}(q_\lambda,{\rm Exp}(c_\lambda),\beta_\lambda)$-bushes as in Theorem \ref{thm1}. Then
$$\frac{1}{\psi^{-1}(\lambda)}Z^{\lambda}_t \rightarrow Z_t\qquad\mbox{in distribution, as $\lambda\rightarrow\infty$, for all $t\ge 0$,}$$
where $(Z_t,t\ge 0)$ is a ${\rm CSBP}(\psi)$ starting from $Z_0=\beta$, with $\psi$ as in {\rm(\ref{offspr})}. If furthermore
$\psi^\prime(0)>-\infty$, then the convergence holds in the almost sure sense.
\end{lemm}
\begin{pf} First consider $t=0$. The initial population sizes are Poisson distributed with
$$\bE(\exp\{-rZ_0^\lambda/\psi^{-1}(\lambda)\})=\exp\{\beta\psi^{-1}(\lambda)(e^{-r/\psi^{-1}(\lambda)}-1)\}\rightarrow e^{-\beta r},$$
as $\lambda\rightarrow\infty$, since $\psi^{-1}(\lambda)\rightarrow\infty$.

For $t>0$, the desired limiting distribution is characterised by (\ref{CSBP and psi solution u}) in terms of $u_t(r)$ and for $x=\beta$. If we integrate (\ref{u solution of psi}), we can identify $u_t(r)$ as the unique solution of
\begin{equation}\label{CSBP u bounds}
    \int^{r}_{u_t(r)}\frac{dv}{{\psi}(v)} = t.
\end{equation}
Consider $Z^{\lambda}_t/{\psi}^{-1}(\lambda)$ as a sum of a Poisson number
of independent ${\rm GW}(q_\lambda,{\rm Exp}(c_\lambda))$-processes. Set
$s=e^{-r/{\psi}^{-1}(\lambda)}$ and apply the branching property at
the first branching time of a single ${\rm GW}(q_\lambda,
{\rm Exp}(c_\lambda))$-tree to obtain for its population size $X_t^\lambda$ at time $t$ with
$w^{\lambda}_t(s) = \bE(s^{X^{\lambda}_t})$
$$%\begin{eqnarray*} 
w^{\lambda}_t(s) \label{w solution}
%    &=& \int^{\infty}_{t} \sum^{\infty}_{0} q_\lambda(k) \mathbf{E}(r^{Z^{\lambda}_t} | V_1=s, %Z^{\lambda}_{V_1}=k )c_\lambda e^{-c_\lambda s} ds  \\
%    & & + \int^{t}_{0} \sum^{\infty}_{0} q_\lambda(k) \mathbf{E}(r^{Z^{\lambda}_t} | V_1=s, %Z^{\lambda}_{V_1}=k )c_\lambda e^{-c_\lambda s} ds \\
    =s e^{-c_\lambda t} + \int_0^t\sum_{k=0}^\infty q_\lambda(k)(\bE(s^{X_{t-y}^\lambda}))^kc_\lambda e^{-c_\lambda y}dy
    =s e^{-c_\lambda t} + \int^{t}_{0} \varphi_{q_\lambda}(w^\lambda_z(s))c_\lambda e^{-c_\lambda (t-z)} dz.
$$%\end{eqnarray*}
Now apply (\ref{offspr}), multiply by $e^{c_\lambda t}$, differentiate with respect to $t$ and rearrange to get
\begin{equation}\label{wlambda}1=\frac{\psi^{-1}(\lambda)\frac{\partial}{\partial t}w_t^\lambda(s)}{\psi(\psi^{-1}(\lambda)(1-w_t^\lambda(s)))}\quad\Rightarrow\quad t=\int_{\psi^{-1}(\lambda)(1-w_t^\lambda(s))}^{\psi^{-1}(\lambda)(1-s)}\frac{dv}{\psi(v)}.\end{equation}
For $s=e^{-r/\psi^{-1}(\lambda)}$, we have $\psi^{-1}(\lambda)(1-s)\rightarrow r$ and by (\ref{CSBP u bounds}) also $\psi^{-1}(\lambda)(1-w_t^\lambda(s))\rightarrow u_t(r)$ and then
$$ \bE(\exp\{-rZ_t^{\lambda}/\psi^{-1}(\lambda)\})=\exp\{-\beta\psi^{-1}(\lambda)(1-w_t^\lambda(s))\}\rightarrow e^{-\beta u_t(r)}\qquad\mbox{as $\lambda\rightarrow\infty$,}
$$
as required.

%-------------------------------------------------------------------------------------------------------------------------------------Proof for A.S.
For the proof of almost sure convergence, recall our notation $\overline{\pi}_t(T)\subset\bU$ for the population alive at time $t$ of a tree $T$, which we will slightly abuse and also apply to bushes. 
%We define sigma-algebras
%$$\mathcal{G}_{\lambda, t}=\sigma\left(\overline{\pi}_s(B_{\lambda^\prime}),\lambda^\prime \geq \lambda, s \leq t \right)\quad\mbox{and}\quad\widetilde{\mathcal{G}}_{\lambda, t}=\sigma\left(\overline{\pi}_s(B_\lambda), s \leq t \right).$$
%It is clear that $\widetilde{\mathcal{G}}_{\lambda, t} \subset \mathcal{\mathcal{G}}_{\lambda, t}$.
Let $0\le\mu<\lambda$, $p=1-\mu/\lambda$ and $\mathcal{G}_{\lambda, t}=\sigma\left\{\overline{\pi}_s(B_{\lambda^\prime}),\lambda^\prime \geq \lambda, s \leq t \right\}$.
Now note that $Z^{\lambda}_t/{\psi}^{-1}(\lambda)$ is $\mathcal{G}_{\lambda, t}$-measurable, and that $\psi^\prime(0) > -\infty$ ensures that $\mathbb{E}(Z^{\mu}_t) < \infty$. Then, a.s.,%Recall $\mathcal{F}^{\mathbb{H}}_t$ in Branching properties of coloured tree in Section 2, fix $\lambda,$
% the branching property is true for all $t>0$ and $\mathcal{F}^{\mathbb{H}, \lambda}$-measurable functions $f_u\colon  \mathbb{T}^\mathbb{H} \rightarrow [0, \infty).$, and for all $\lambda \in [0, \infty)$, the branching property is true adapted by $\mathcal{F}^{\mathbb{H}, \lambda}$ respectively. Moreover for fixed $\lambda$ and $t,$ it is easily understood (9) still hold, and on the other hand $\tilde{\mathcal{G}}_{\lambda, t} \subset \mathcal{F}^{\mathcal{H}, \lambda}.$ And according to (9),
$$%\begin{eqnarray*}
\mathbb{E}\left(\left.Z^{\mu}_{t} \right| \mathcal{G}_{\lambda, t} \right)
%&=&\mathbb{E}\left(\sum_{u \in B_\lambda}\mathbf{1}_{\{u \in \overline{\pi}_t(B_\lambda),
%\gamma_u(B_\lambda^{p-\rm col})=0\}} | \mathcal{G}_{\lambda, t} \right)\\
=\mathbb{E}\left(\left.\sum_{u \in
\overline{\pi}_t(B_\lambda)}\mathbf{1}_{\{\gamma_u(B_\lambda^{p-\rm col})=0\}} \right|
\mathcal{G}_{\lambda, t} \right)
%=\mathbb{E}\left(\left.\sum_{u \in
%\overline{\pi}^\lambda_t(B_\lambda)}\mathbf{1}_{\{\gamma_u(B_\lambda^{p-\rm col}) =0\}} \right|
%\widetilde{\mathcal{G}}_{\lambda, t} \right),
= \sum_{u \in
\overline{\pi}_t(B_\lambda)}\mathbb{E}\left(\left.\mathbf{1}_{\{\gamma_u(B_\lambda^{p-\rm col})=0\}}\right|\mathcal{G}_{\lambda, t} \right),
$$%\end{eqnarray*}
since $Z_t^\mu=\#\overline{\pi}_t(B_\mu)$. By the branching property in Remark \ref{Markov}(c), applied to functions
$f_u(\cdot)=\mathbf{1}_{\{\gamma_\emptyset(\cdot)=0\}},$ and $f_v(\cdot)\equiv 1$ for $v \neq u,$ we obtain, a.s.,
%\begin{eqnarray*}
$$
%\mathbb{E}\left(\left.\sum_{u \in
%\overline{\pi}_t(B_\lambda)}\mathbf{1}_{\{\gamma_u(B_\lambda^{p-\rm col}) =0\}} \right|
%\widetilde{\mathcal{G}}_{\lambda, t} \right)&=& 
\sum_{u \in
\overline{\pi}_t(B_\lambda)}\mathbb{E}\left(\left.\mathbf{1}_{\{\gamma_u(B_\lambda^{p-\rm col})=0\}}\right| \mathcal{G}_{\lambda, t} \right)
=Z^{\lambda,\beta}_t \mathbb{P}\left(\gamma_\emptyset(B_\lambda^{p-\rm col})=0\right)
=(1-g_\lambda(p))Z^{\lambda,\beta}_t.
$$%\end{eqnarray*}
According to (\ref{g in psi}), this shows that
$$\mathbb{E}\left(\left.\frac{Z^{\mu}_t}{\psi^{-1}(\mu)} \right| \mathcal{G}_{\lambda, t} \right) = \frac{Z^{\lambda}_t}{{\psi}^{-1}(\lambda)}\qquad\mbox{a.s.,}$$
which is the martingale property in the (decreasing) filtration $(\mathcal{G}_{\lambda,t},\lambda\ge 0)$ that implies that $Z^{\lambda}_{t}/\psi^{-1}(\lambda) \rightarrow Z_t$ almost surely as $\lambda \rightarrow
\infty$, see e.g. \cite{Wil-mart}.
\end{pf}

%Actually, the assumption $\psi^\prime(0+)>-\infty$ can be dropped and even the explosive case can be included if we interpret the statement as convergence before the explosion time $\dagger=\inf\{t\ge 0:Z_t=\infty\}$ and divergence to $\infty$ after $\dagger$. Specifically, to a 
%branching mechanism $\psi$ with $\psi^\prime(0+)=-\infty$ and $h>0$ associate a truncated branching mechanism $\psi_{|(0,h)}$ such 
%that 
%$$\psi_{|(0,h)}(r)=br+ar^2+\int_{(0,h)}(e^{-rx}-1+rx{\bf 1}_{\{x<1\}})\Pi(dx).$$
%It is not hard to show, e.g. via coupling of L\'evy processes and Lamperti's time-change \cite{Lam-CSBP}, that we can couple families 
%$(Z^{\lambda,\beta,h}_t,t\ge 0)$, $h>0$, $\lambda\ge 0$, such that $\bP(Z^{\lambda,\beta,h}_t=Z^{\lambda,\beta}_t\mbox{ for all $\lambda\ge 0$})\rightarrow\bP(\dagger>t)$, as $h\rightarrow\infty$. Similarly, by considering the ladder time processes of L\'evy processes \cite[Chapter VII]{Ber-LP}, we obtain that $(\psi^{-1}_{|(0,h)}(\lambda)-\psi^{-1}(\lambda),\lambda\ge 0)$ is bounded for each $h>0$, and so $\psi^{-1}(\lambda)/\psi_{|(0,h)}^{-1}(\lambda)\rightarrow 1$ as $\lambda\rightarrow\infty$. We leave the details to the reader.

%------------------------------------------------------------------------------------------------------------------Consistent kappa bush

\section{Growth of ${\rm GW}(q_\lambda, \kappa_\lambda, \beta_\lambda)$-bushes: lifetimes}\label{Consistend kappa bush}

Theorem \ref{thm2} considers Galton-Watson bushes $B_\lambda$ with general
independent $\kappa_\lambda$-distributed lifetimes. The main statement beyond the exponential case of Theorem \ref{thm1} is that consistency of a full
family $(B_\lambda,\lambda\ge 0)$ under Bernoulli leaf colouring
requires $\kappa_\lambda$ to be geometrically divisible.
%--------------------------------------------------------------------------------------------------------------------------Growth for general edge trees
\subsection{Proof of Theorem \ref{thm2}}\label{Growth for general edge trees}
${\rm (i)}\Rightarrow{\rm (ii)}$: Suppose, (i) holds. In particular, $q_\mu$ is then the $(1-\mu/\lambda)$-reduced
offspring distribution associated with $q_\lambda$, for all $\mu<\lambda$. By Theorem \ref{thm1},
$\varphi_q(s)=s+\widetilde{\psi}(1-s)$, where $\widetilde{\psi}$ has the form (\ref{lkbm}). To be specific, let $c=1$,
$\beta\in(0,\infty)$ and parametrise $(q_\lambda,\beta_\lambda)$ using $\psi$ as in Section \ref{Consistent family}.

Consider a ${\rm GW}(q_\lambda,\kappa_\lambda)$-tree $T_\lambda$. By Remark \ref{alphakappa}, Lemma \ref{alpha and psi} applies. Using notation
similar to the reduction procedure of Section \ref{coloured leaves and 2-colour BP} for $T=T_\lambda$, $p=1-1/\lambda$, $\iota=\iota_\lambda$, $\lambda\ge 1$,
we can write $\zeta_\emptyset\sim \kappa_1$ as
\begin{equation}\label{writeas}\zeta_\emptyset=\sum_{v\in \iota_\lambda^{-1}(\emptyset)}\zeta_v^\lambda\overset{{\rm (d)}}{=} \sum_{j=1}^{G^{(\alpha)}}X_\alpha^{(j)},
\end{equation}
for $X_\alpha^{(j)}\sim\kappa_\lambda$, $j\ge 1$, and $G^{(\alpha)} \sim {\rm Geo}(\alpha)$ with
$\alpha=\psi^\prime(\psi^{-1}(1))/\psi^\prime(\psi^{-1}(\lambda))$ independent. Specifically, consider $\Gamma_n=\{G_\emptyset^{(\alpha)}=n\}$, where $G_\emptyset^{(\alpha)}=\iota^{-1}_\lambda(\emptyset)$ in the notation of Section \ref{Recon for GW expo bushes1}. On $\Gamma_n$, write $v_0=\emptyset$ and for
$1\le j\le n-1$, denote by $v_j$ the unique black child of $v_{j-1}$, then an $n$-fold inductive application of
the two-colours branching property of Section \ref{coloured leaves and 2-colour BP}, also summing over all offspring numbers and colour combinations
as in Lemma \ref{alpha and psi}, yields for $p=1-1/\lambda$
\beq \bE_{q_\lambda\otimes\kappa_\lambda}^{p-\rm col}\!\!\left[\left.\prod_{j=0}^{n-1} k_j(\zeta_{v_j});\Gamma_n\right|\gamma_\emptyset\!=\!0\right]
        &\!\!\!\!\!=\!\!\!\!\!&\int_{(0,\infty)}\!\!\!\!k_0(z)\kappa_\lambda(dz)\varphi_{q_\lambda}^\prime(g_\lambda(p))\bE_{q_\lambda\otimes\kappa_\lambda}^{p-\rm col}\!\!\left[\left.\prod_{j=0}^{n-2}k_{j+1}(\zeta_{v_{j}});\Gamma_{n-1}\right|\gamma_\emptyset\!=\!0\right]\\
        &\!\!\!\!\!=\!\!\!\!\!&\left(\prod_{j=0}^{n-1}\int_{(0,\infty)}k_j(z)\kappa_\lambda(dz)\right)(1-\alpha)^{n-1}\alpha.
\eeq
By Corollary \ref{9},  $\alpha\downarrow\psi^\prime(\psi^{-1}(1))/\psi^\prime(\infty)=1/\widetilde{\psi}^\prime(\infty)$ as $\lambda\rightarrow\infty$.
Therefore, we can write $\zeta_\emptyset\sim\kappa_1$ as in (\ref{writeas}) for all $\alpha>1/\widetilde{\psi}^\prime(\infty)$, with the convention $1/\widetilde{\psi}^\prime(0)=0$ if $\widetilde{\psi}^\prime(0)=\infty$. This yields (ii).

${\rm (ii)}\Rightarrow{\rm (i)}$: Suppose, (ii) holds. By Theorem \ref{thm1}, the family $(q_\lambda,\lambda\ge 0)$
exists as required. By Theorem \ref{thm1} and Remark \ref{alphakappa}, we can express $\varphi_{q_\lambda}$ and
$\alpha$ in terms of $\psi$ as in Section \ref{Consistent family} and Lemma \ref{alpha and psi}, choosing $c=1$.

For $\lambda>1$, geometric divisibility of $\kappa$ permits us to define $\kappa_\lambda$ as the
distribution of
$X_\alpha^{(j)}$
for $\alpha=\psi^\prime(\psi^{-1}(1))/\psi^\prime(\psi^{-1}(\lambda))=1/\widetilde{\psi}^\prime(\widetilde{\psi}^{-1}(\lambda\widetilde{\psi}(1)))>1/\widetilde{\psi}^\prime(\infty)$.
Now consider a ${\rm GW}(q_\lambda,\kappa_\lambda)$-tree and $p=1-1/\lambda$. Use notation from ${\rm (i)}\Rightarrow{\rm (ii)}$ and also set
$v^*=v_{G^{(\alpha)}_\emptyset-1}$. On $\Upsilon_j=\{\nu_{v^*}-\gamma_{v^*1}-\cdots-\gamma_{v^*\nu_{v^*}}=j\}$ for $j\ge 2$,
denote by $w_1,\ldots,w_j$ the black children of $v^*$. Then, by Remark \ref{alphakappa} and repeated application of the two-colours branching property,
we obtain
\beq &&\hspace{-0.5cm}\bE_{q_\lambda\otimes\kappa_\lambda}^{p-\rm col}\left[\left.\exp\left\{-r\sum_{m=0}^{G_\emptyset^{(\alpha)}-1}\zeta_{v_m}\right\}\prod_{i=1}^jf_i(\overline{\theta}_{w_i});\Upsilon_j\right|\gamma_\emptyset=0\right]\\
      &&=\sum_{n=1}^\infty\left(\int_{(0,\infty)}e^{-rz}\kappa_\lambda(dz)\right)^n\ (1-\alpha)^{n-1}\ \varphi_{q_\lambda}^{(j)}(g_\lambda(p))\frac{(1-g_\lambda(p))^{j-1}}{j!}\ \prod_{i=1}^j\bE_{q_\lambda\otimes\kappa_\lambda}^{p-\rm col}\left[\left.f_i\right|\gamma_\emptyset=0\right].
\eeq This is the branching property characterizing ${\rm
GW}(q,\kappa)$, because the first term is the Laplace transform of a
geometric sum with distribution $\kappa$, up to a factor of
$\alpha$, and for the middle term we identify the offspring
distribution $q$ %(cf. \cite{DuW}) 
using all the cancellations due to
(\ref{offspr}), (\ref{g in psi}) and (\ref{imporesult})
\beq &&\hspace{-0.7cm}\sum_{j=2}^\infty s^j\frac{1}{\alpha}\frac{\varphi_{q_\lambda}^{(j)}(g_\lambda(p))}{j!}(1-g_\lambda(p))^{j-1}+\left(1-\sum_{j=2}^\infty\frac{1}{\alpha}\frac{\varphi_{q_\lambda}^{(j)}(g_\lambda(p))}{j!}(1-g_\lambda(p))^{j-1}\right)\\
    &&\hspace{-0.2cm}=1+\frac{\varphi_{q_{\lambda}}(g_\lambda(p)+s(1-g_\lambda(p)))-1-(1-s)\varphi_{q_\lambda}^\prime(g_\lambda(p))(1-g_\lambda(p))}{\alpha(1-g_\lambda(p))}=s+\frac{\psi(\psi^{-1}(1)(1-s))}{\psi^{-1}(1)\psi^\prime(\psi^{-1}(1))}
\eeq
confirming that $(q,\kappa)$ is the $(1-1/\lambda)$-reduced pair associated with $(q_\lambda,\kappa_\lambda)$. For $\mu<1$, set $\alpha=\psi^\prime(\psi^{-1}(\mu))/\psi^\prime(\psi^{-1}(1))$ and define $\kappa_\mu$
to be the distribution of
$$\sum_{j=1}^{G^{(\alpha)}}X^{(j)},\qquad\mbox{for independent $X^{(j)}\sim\kappa$, $j\ge 1$, independent of $G^{(\alpha)}\sim{\rm Geo}(\alpha)$.}$$
As above, $(q_\mu,\kappa_\mu)$ is the
$(1-\mu)$-reduced pair associated with $(q,\kappa)$. The reduction relation for $0\le\mu<\lambda<\infty$ follows
using transitivity of colouring reduction (see Remark \ref{remdet}(a)) for $0\le\nu<\mu<\lambda<\infty$. Specifically, for $\mu=1$, this yields that $(B_\lambda^{(1-\nu/\lambda)-{\rm rdc}},B_\lambda)\overset{{\rm (d)}}{=}  (B_\nu,B_\lambda)$. For $\nu=1<\mu<\lambda$ and $\nu<\mu<\lambda=1$, this argument can be combined with the uniqueness of the divisor distribution (Lemma \ref{uniqgeodiv}). This completes the proof of ${\rm (ii)}\Rightarrow{\rm (i)}$.

We identified $(\beta_\lambda,\lambda\ge 0)$ in ${\rm(i)}\Rightarrow{\rm (ii)}$. The same reasoning as in Remark \ref{alphakappa} allows us to combine (i) here and Theorem \ref{thm1} to see that $(q,\kappa,\beta)$ is the $(1-1/\lambda)$-reduced triplet associated with $(q_\lambda,\kappa_\lambda,\beta_\lambda)$. The existence of $(B_\lambda,\lambda\ge 0)$ now follows from Kolmogorov's consistency theorem. Uniqueness of the families $(q_\lambda,\lambda\ge 0)$, $(\kappa_\lambda,\lambda\ge 0)$ and $(\beta_\lambda,\lambda\ge 0)$ for each
$\beta=\beta_1\in(0,\infty)$ follows from the uniqueness results in Theorem \ref{thm1} and as shown in ${\rm (ii)}\Rightarrow{\rm (i)}$.\hspace*{\fill} $\square$

% \begin{flushright}$\Box$\end{flushright}

%--------------------------------------------------------------------------------------------------------------------------------------Reconstruction General
\subsection{Reconstruction procedures and backbone decomposition}\label{Reconstruction General Edge}
If $\kappa_\lambda(dz)=f_\lambda(z)dz$ is absolutely continuous for all $\lambda\ge 0$ and $\zeta=\zeta_1+\cdots+\zeta_{G^{(\alpha)}}\sim\kappa_\mu$ for independent $\zeta_j\sim\kappa_\lambda$ and $G^{(\alpha)}\sim{\rm Geo}(\alpha)$ for $\alpha$ as in Section \ref{Recon for GW expo bushes1}, we find conditional joint
distributions as in (\ref{joint exponential})
%. Since $\zeta_{u_j},
%j=1,2,\ldots $ are independent identical distributed as $\zeta_{u_1}
%\sim \kappa_\lambda$, the conditional joint distribution can be
%expressed as follows:
\begin{equation}\label{joint general edge}
f_{\zeta_{1},\ldots,\zeta_{n-1}|G^{(\alpha)}=n,\zeta=z}(y_1,\ldots,y_{n-1})=
\frac{f_\lambda(z-\sum^{n-1}_{j=1}y_j)\prod^{n-1}_{j=1}f_\lambda(y_j)}{f_\lambda^{*(n)}(z)}
\end{equation}
for $y_j>0$ with $y_1+\cdots+y_{n-1}<z$, $n\ge 1$, where
$f_\lambda^{*(n)}$ is the $n$th convolution power of $f_\lambda$.

\subsubsection{Reconstruction procedure for ${\rm GW}(q_\lambda, \kappa_\lambda,\beta_\lambda)$-bushes}

For $B_\mu\sim{\rm GW}(q_\mu,\kappa_\mu,\beta_\mu)$ the
procedure in Section \ref{Recon for GW expo bushes2} with ${\rm Exp}(c_\lambda)$ replaced by $\kappa_\lambda$ and
(\ref{joint exponential}) replaced by (\ref{joint general edge}) constructs $B_\lambda\sim{\rm GW}(q_\lambda,\kappa_\lambda,\beta_\lambda)$.

In the general case, one could use regular conditional distributions. Alternatively, we can adapt the procedure in
Section \ref{Recon for GW expo bushes2} using the subordinator (random walk) representation of geometrically infinitely
(finitely) divisible distributions as explained below.

\subsubsection{Reconstruction procedure with subordinators in the case where $\kappa$ is g.i.d.}\label{recprocgid}

In the g.i.d.\ case, let $\bH=\bigcup_{\zeta\in(0,\infty)}\{\zeta\}\times\bD([0,\zeta],[0,\infty))$,
where $\bD([0,\zeta],[0,\infty))$ is the set of functions $f\colon[0,\zeta]\rightarrow[0,\infty)$ that are right-continuous with left limits, equipped with the Borel sigma-algebra generated by the metric topology induced by $$d((\zeta_1,f_1),(\zeta_2,f_2))=|\zeta_1-\zeta_2|+d_{\rm Sk}(f_1(\cdot\wedge\zeta_1),f_2(\cdot\wedge\zeta_2)),$$ where $d_{\rm Sk}$ is a
metric that generates Skorohod's topology on $\bD([0,\infty),[0,\infty))$, see e.g. \cite[Section IV.1]{LeG-snake}. Consider a subordinator $(\sigma(t),t\ge 0)$, under $\bP$, such that $\sigma(V)\sim\kappa_1$ for an
independent $V\sim{\rm Exp}(1)$. Define the following measure on $\bN_0\times\bH$
\begin{equation}\label{Qmu}\bQ_\mu(\{j\}\times H)=q_\mu(j)\int_{(0,\infty)}\bP((z,(\sigma(t),0\le t\le z))\in H)c_\mu e^{-c_\mu z}dz.\end{equation}
Now consider a bush $\overline{B}_\mu$ of $N_\mu\sim{\rm Poi}(\beta_\mu)$ random trees with distribution $\bP_{\bQ_\mu}$ as defined in the branching property of Section \ref{Lifetime marks and discr BP in cont. time} for the measure $\bQ_\mu$ just defined.
Then the reconstruction procedure of Section \ref{Recon for GW expo bushes2} can be applied subdividing subordinator
lifetimes $\zeta_u^{(i)}\sim{\rm Exp}(c_\mu)$ rather than directly the population lifetimes $\sigma_u^{(i)}(\zeta_u^{(i)})\sim\kappa_\mu$. Also define
$$\sigma_{u,m}^{(i)}(t)=\sigma_u^{(i)}(\zeta_{u,1}^{(i)}+\cdots+\zeta_{u,m-1}^{(i)}+t)-\sigma_u^{(i)}(\zeta_{u,1}^{(i)}+\cdots+\zeta_{u,m-1}^{(i)}),\quad 0\le t\le\zeta_{u,m}^{(i)}, 1\le m\le G_u^{(\alpha,i)}.$$
The remainder of the procedure is easily adapted. The resulting $\overline{B}_\lambda$ is a bush of $N_\lambda\sim{\rm Poi}(\beta_\lambda)$ random trees with distribution $\bP_{\bQ_\lambda}$.

\subsubsection{Reconstruction procedure with random walks in the case where $\psi^\prime(\infty)<\infty$.}

In the case where $\kappa$ is geometrically divisible up to $\alpha_0(1)=\psi^\prime(\psi^{-1}(1))/\psi^\prime(\infty)>0$, let $\bH=\bigcup_{n\in\mathbb{N}}\{n\}\times[0,\infty)^{n+1}$ be the space of random walk paths. On $\bN_0\times\bH$, consider
$$\bQ^{\rm RW}_\mu(\{j\}\times H)=q_\mu(j)\sum_{n=1}^\infty\bP((n,(\sigma(k),0\le k\le n))\in H)(1-\alpha_0(\mu))^{n-1}\alpha_0(\mu),$$
where, under $\bP$, $(\sigma(k),k\ge 0)$ is a random walk with $\sigma(k+1)-\sigma(k)\sim\kappa_\infty$ i.i.d.
In fact, in this case, the distribution $\kappa_\infty$ exists since the Laplace transforms of $\kappa_\lambda$
converge as $\lambda\rightarrow\infty$ to a completely monotone function continuous at zero. Then the reconstruction
procedure of Section \ref{Recon for GW expo bushes2} can be applied subdividing geometric ``random walk lifetimes'' $G_u^{(i)}\sim{\rm Geo}(\alpha_0(\mu))$ into
$G_u^{(i)}=G_{u,1}^{(i)}+\cdots+G_{u,G_u^{(\alpha,i)}}$, for independent $G_{u,m}^{(i)}\sim{\rm Geo}(\alpha_0(\lambda))$, $m\ge 1$, and $G_u^{(\alpha,i)}\sim{\rm Geo}(\alpha)$.
Also define
$$\sigma_{u,m}^{(i)}(k)=\sigma_u^{(i)}(G_{u,1}^{(i)}+\cdots+G_{u,m-1}^{(i)}+k)-\sigma_u^{(i)}(G_{u,1}^{(i)}+\cdots+G_{u,m-1}^{(i)}),\qquad 0\le k\le G_{u,m}^{(i)}.$$
The remainder of the procedure is easily adapted. The resulting $\overline{B}_\lambda$ is a bush of $N_\lambda\sim{\rm Poi}(\beta_\lambda)$ random trees with distribution $\bP_{\bQ^{\rm RW}_\lambda}$.

\subsubsection{Backbone decomposition of supercritical Bellman-Harris processes}

The reconstruction procedures that build a ${\rm GW}(q_\lambda,\kappa_\lambda,\beta_\lambda)$-bush $B_\lambda$ from a ${\rm GW}(q_\mu,\kappa_\mu,\beta_\mu)$-bush give rise to decompositions of the associated Bellman-Harris process $Z_t^\lambda=\#\overline{\pi}_t(B_\lambda)$ along the ${\rm GW}(q_\mu,\kappa_\mu,\beta_\mu)$-bush $B_\mu$. In the sequel, we will write ${\rm BH}(q_\lambda,\kappa_\lambda)$ for such a Bellman-Harris process and when we specify its initial distribution, all these individuals are taken with zero age.

In the supercritical case $\psi^\prime(0)<0$, note that the ${\rm GW}(q_{\lambda,\rm red}^{(1-\mu/\lambda)-\rm rdc},\kappa_\lambda)$-trees
with offspring distribution as in (\ref{red offspring distribution}) that are grafted onto $B_\mu$ are subcritical for all $0\le\mu<\lambda$. The case $\mu=0$ in is at the heart of many decompositions in various settings, mainly continuous analogues with and without spatial motion, see \cite{BKM-10,DuW,EtW-99,EvO-91}. As an immediate consequence of our reconstruction procedures, we obtain a version of the backbone decomposition for Bellman-Harris processes. 	

\begin{coro}\label{cor14} Let $\psi$ be a supercritical %(non-explosive) 
branching mechanism, $B_0$ a bush of $N_0\sim{\rm Poi}(\beta\psi^{-1}(0))$ random trees with distribution 
  $\bP_{\bQ_0}$ as in {\rm (\ref{Qmu})}. Subdivide each subordinator lifetime as in Section \ref{recprocgid} to get a bush $\widehat{B}_0$. For each $u\in\widehat{B}_0$
  independently, given $\widehat{\nu}_u=m$ children, consider a ${\rm BH}(q_{\lambda,\rm red}^{1-\rm rdc},\kappa_\lambda)$-process $Z^{(u)}$ with $Z^{(u)}_0$ 
  of distribution {\rm(\ref{redgraft})}. Also consider a ${\rm BH}(q_{\lambda,\rm red}^{1-\rm rdc},\kappa_\lambda)$-process $Z^{\rm root}$ with $Z^{\rm root}_0\sim{\rm Poi}(\beta(\psi^{-1}(\lambda)-\psi^{-1}(0)))$. Then the process
	$$Z_t=\#\overline{\pi}_t(B_0)+Z^{\rm root}_t+\sum_{u\in\widehat{B}_0\colon \omega_u\le t}Z^{(u)}_{t-\omega_u}\vspace{-0.2cm}$$
  is a ${\rm BH}(q_\lambda,\kappa_\lambda)$-process with $Z_0\sim{\rm Poi}(\beta\psi^{-1}(\lambda))$.
\end{coro}

\subsection{Limiting trees and branching processes as $\lambda\rightarrow\infty$}\label{realtree}

\subsubsection{Convergence of trees: $\bR$-tree representations, L\'evy trees and snakes}\label{secsnake}

A random marked tree $\overline{T}_\lambda$ with distribution $\bP_{\bQ_\lambda}$ as in Section \ref{recprocgid} specifies marks $\zeta_u\sim{\rm Exp}(c_\lambda)$ as well as $\sigma_u(\zeta_u)\sim\kappa_\lambda$. Therefore, we can associate coupled
trees
\begin{eqnarray*}&&T_\lambda^\circ=\{(u,\zeta_u)\colon (u,\zeta_u,\sigma_u)\in\overline{T}_\lambda\}\sim{\rm GW}(q_\lambda,{\rm Exp}(c_\lambda),\beta_\lambda)\\
    \mbox{and}&&T_\lambda^\bullet=\{(u,\sigma_u(\zeta_u))\colon (u,\zeta_u,\sigma_u)\in\overline{T}_\lambda\}\sim{\rm GW}(q_\lambda,\kappa_\lambda,\beta_\lambda),
\end{eqnarray*}
that only differ in their lifetimes. In the same way, we obtain coupled bushes and consistent families $(B_\lambda^\circ,\lambda\ge 0)$ and $(B_\lambda^\bullet,\lambda\ge 0)$. Several other representations of $(T_\lambda^\circ,T_\lambda^\bullet)$ are natural. For an $\bR$-tree representation
$$\cT_\lambda^\circ=\{\rho\}\cup\bigcup_{u\in\tau_\lambda}\{u\}\times(\alpha_u^\circ,\omega_u^\circ]$$
of $T_\lambda^\circ$, with root $\rho=(\emptyset,0)$ and metric $\mathrm{d}$ given by $\mathrm{d}((v,s),(w,t))=|t-s|$ for $v\preceq w$ or $w\preceq v$,
$$\mathrm{d}((uiv,s),(ujw,t))=s+t-2\omega_u^\circ\qquad\mbox{for $u,v,w\in\bU$ and $i,j\in\bN$, $i\neq j$,}$$
see \cite[Sect. 3.3]{DuW}, we can consider the measure
$$\cW_\lambda(\{u\}\times(a,b])=\sigma_u(b-\alpha_u^\circ)-\sigma_u(a-\alpha_u^\circ),\qquad\alpha_u^\circ\le a<b\le\omega_u^\circ,u\in\tau_\lambda,$$
which for $\sigma_u\overset{{\rm (d)}}{=} \sigma$ as in (\ref{LaExponent})-(\ref{LK1}), $u\in\tau_\lambda$, is of the form $\cW_\lambda=d{\rm Leb}+\cR_\lambda$, where $\cR_\lambda$ is an infinitely divisible independently scattered random measure on $\cT_\lambda^\circ$, in the sense of \cite{Kal,RaR,WoU}, with $F$-measure ${\rm Leb}\otimes\Lambda$, where ${\rm Leb}$ denotes Lebesgue measure on $\cT_\lambda^\circ$. 
For technical convenience, we can embed $(\cT_\lambda^\circ,\mathrm{d},\rho)$ as a metric subspace of  $\ell_1(\bN)=\{x=(x_n)_{n\ge 1}\colon x_n\in[0,\infty),n\ge 1;||x||_1<\infty\}$, where $||x||_1=\sum|x_n|$ is the $\ell_1$-norm,
with root $0\in\ell_1(\bN)$. We can embed bushes $(\cB_\lambda^\circ,\lambda\ge 0)$ consistently following \cite[Remarks 4.9-4.10]{DuW}: in the (sub)critical case, this family starts from $\varnothing$, is piecewise constant and evolves by adding single branches, so we can represent the $j$th branch, of length $L_j^\circ$ say, as $[[x(j),x(j)+L_j^\circ e_j]]$ for some $x(j)=(x_1(j),\ldots,x_{j-1}(j),0,\ldots)\in\cB_{\lambda_{j-1}}^\circ$ and $e_j=(0,\ldots,0,1,0,\ldots)$ the $j$th coordinate vector in $\ell_1(\bN)$. Similarly, we can embed $(\cB_\lambda^\bullet,\lambda\ge 0)$ for certain $L_j^\bullet=\sigma_j(L_j^\circ)$. The supercritical case can be handled 
using \cite[Proposition 3.7]{DuW}.

For the $\ell_1(\bN)$-embedding, we can consider $\cW_\lambda$ as a random measure on $\ell_1(\bN)$ (with support included in the embedded $\cB_\lambda^\circ$).
Then $\cB_\lambda^{\circ,\rm wt}=(\cB_\lambda^\circ,\cW_\lambda)$ is a weighted $\bR$-tree in the sense of \cite{EvW-07,GPW-09}, in the (sub)critical case, but weak or vague convergence as $\lambda\rightarrow\infty$ are not appropriate since in the limit the measures become infinite on any ball around a point in some $\cB_\lambda^\circ$. Nevertheless, consistency implies that in the case where convergence to a locally compact and separable L\'evy bush $\cB^\circ$ occurs in the Gromov-Hausdorff sense \cite[Theorem 5.1]{DuW}, the random measures $\cW_\lambda$ consistently build an infinitely divisible independently scattered random measure $\cW$ on $\cB^\circ\subset\ell_1(\bN)$ whose Poissonian component still has $F$-measure ${\rm Leb}|_{\cB^\circ}\otimes\Lambda$, while the continuous component is still $d{\rm Leb}|_{\cB^\circ}$ and where ${\rm Leb}$ is one-dimensional Lebesgue measure on $\ell_1(\bN)$. 

For $x\in\cB^\circ$ consider the path $[[0,x]]=\{y\in\cB^\circ\colon y_n\le x_n,n\ge 1\}=\{f_x(t)\colon 0\le t\le||x||_1\}$, where  $||f_x(t)||_1=t$, $0\le t\le||x||_1$. Then $\cS\colon\cB^\circ\rightarrow\bH$ with $\bH$ as in Section \ref{recprocgid} given by $\cS(x)=(||x||_1,(\cW([[0,f_x(t)]]),0\le t\le||x||_1))$ can be seen as a snake in the sense of \cite[Section I.3.2]{LeG-snake}, where the
spatial motion here is a subordinator with characteristic pair $(d,\Lambda)$ as in (\ref{LK1}). See \cite{JaM} for another setting where discontinuous snakes appear naturally.

\begin{prop} In the setting of Theorem \ref{thm2}, if $\cB^\circ$ is separable and locally compact, then 
  $\cB^\bullet$, the closure of $\bigcup\cB^\bullet_\lambda$ in $\ell_1(\bN)$, is separable, but locally compact only if $\Lambda=0$. Also, if $\cB^\circ$ is furthermore bounded, then $\cB^\bullet$ is bounded if and only if $\Lambda=0$.
\end{prop}
\begin{pf} Separability of $\cB^\bullet$ is trivial as $\cB^\bullet_\lambda$ is separable. Now argue conditionally given $(\cB^\circ_\lambda)_{\lambda\ge 0}$. 
 To show that boundedness and local compactness fail unless $\Lambda=0$, first consider $\Lambda=\delta_h$, the Dirac measure in $h>0$. It is not hard to
  show that the subtree of $\cB^\circ$ above every vertex, except leaves, has infinite total length. As a consequence, there will be
  $\lambda_1>0$ for which $\cW_{\lambda_1}$ has an atom $x_1$, and the subtree $\{y\in\cB^\circ\colon x_1\in[[0,y[[\}$ above $x_1$ has 
  infinite length. Inductively, we find an increasing sequence $x_1\prec x_2\prec\cdots$ of atoms of $\cW$ showing that
  $\cB^\bullet$ is not bounded. 
  
  Now assume that $\cB^\bullet$ is locally compact. By the Hopf-Rinow theorem \cite{Gro}, closed balls are then compact. Since $\cB^\circ$ is locally compact, \cite[Remark 5.1]{DuW} implies that \cite[Theorem 5.1]{DuW} applies. In particular, $\psi^\prime(\infty)=\infty$, so that $\cB^\circ\setminus[[0,L^\circ_1e_1]]$ has infinitely many connected
  components, each of infinite total length, so that they all contain atoms of $\cW$. However, $L_1^\bullet=\sigma_1(L_1^\circ)$ is bounded and for a large enough ball with $B(r)=\{x\in\cB^\bullet\colon ||x||_1\le r\}\supset[[0,L^\bullet_1e_1]]$, the set $B(r+h)\setminus[[0,L^\bullet_1e_1]]$ also has infinitely many connected components each exceeding diameter $h$. Considering a cover of $B(r+h)$ by open balls of radii less than $h/2$, there is no finite subcover, as each connected component needs at least one ball that does not intersect any other connected component, which contradicts the compactness of $B(r+h)$. So $\cB^\bullet$ is not locally compact.

  For any $\Lambda\neq 0$ and $h>0$ with $\Lambda((h,\infty))>0$, let 
  $\widetilde{\Lambda}=\Lambda((h,\infty))\delta_h$, couple $\cB^\bullet$ and 
  $\widetilde{\cB}^\bullet$ and argue as above that $\widetilde{\cB}^\bullet$ is not bounded nor locally 
  compact, then deduce the same for $\cB^\bullet$. 
\end{pf}

Since $\cB^\bullet$ is separable for a large class of branching mechanisms and lifetime subordinators, the framework of
\cite{GPW-09} can be used to further study these trees.

\subsubsection{Superprocesses, backbones and convergence of Bellman-Harris processes}\label{CSBPKpsi}

Sagitov \cite{Sag-91} studied convergence of Bellman-Harris processes to certain non-Markovian CSBP whose distribution is best described via Markovian superprocesses that record residual lifetimes, see also \cite{KaS-98}. Specifically, a $(\xi,K,\psi)$-superprocess \cite{Dyn-91} is a Markov process $M=(M_t,t\ge 0)$ on the space of finite Borel measures $\bM([0,\infty))$ with transition semigroup characterised by
\begin{equation}\label{meassg}\bE\left(\left.\!\exp\left\{\!-\!\!\int_{[0,\infty)}\!\!f(z)M_t(dz)\!\right\}\right|M_0\!=\!m\!\right)\!=\exp\left\{\!-\!\!\int_{[0,t]}\!\!u_{t-z}(f)m(dz)\!-\!\!\int_{(t,\infty)}\!\!f(z\!-\!t)m(dz)\!\right\},\hspace{-0.0cm}\end{equation}%\pagebreak
for all $f\colon[0,\infty)\rightarrow[0,\infty)$ bounded continuous, where $u_t(f)$ is the unique (at least if $\psi^\prime(0)>-\infty$) nonnegative solution of 
\begin{equation}\label{measint}u_t(f)+\int_{[0,t]}\psi(u_{t-s}(f))dH_s=\bE(f(\xi_t))\end{equation}
and $H_s=\bE(K_s)$ is the renewal function of a strictly increasing subordinator or random walk $\sigma$, %of Section \ref{recprocgid}
here in terms of the local time process $K_s=\inf\{t\ge 0\colon \sigma(t)>s\}$, which is an additive functional of our particular choice $\xi_t=\sigma(K_t)-t$ of Markovian spatial motion. Since branching only occurs at the origin, $M$ is a catalytic superprocess \cite{DaF-94,Dyn-95}. We call the associated population size process $Z_t=M_t([0,\infty))$ a ${\rm CSBP}(K,\psi)$. Such processes appear as limits of Bellman-Harris processes, also in our setting. Unless stated otherwise, we understand $Z_0=\beta$ as $M_0=\beta\delta_0$. 

\begin{prop}\label{Sagconv} Let $Z_t^{\lambda}=\#\overline{\pi}_t(B_\lambda)$ for a consistent family $(B_\lambda)_{\lambda\ge 0}$ of
${\rm GW}(q_\lambda,\kappa_\lambda,\beta_\lambda)$-bushes with branching mechanism $\psi$ as in Theorem \ref{thm2}. Suppose that $\psi$ is subcritical or critical, i.e. $\psi^\prime(0)\ge 0$, and that $\psi^\prime(\infty)=\infty$. Let $\sigma$ be a subordinator as in Section \ref{recprocgid}. Then, for all $t \geq 0$
$$\frac{1}{\psi^{-1}(\lambda)}Z^{\lambda}_t \longrightarrow Z_t\qquad\mbox{almost surely, as $\lambda\rightarrow\infty$,}$$
where $(Z_t,t\ge 0)$ is a ${\rm CSBP}(K,\psi)$ starting from $Z_0=\beta$.
\end{prop}
\begin{pf} The almost sure convergence follows from the same martingale argument as in Lemma \ref{CSBP Limit Lemma x}. The identification of the limiting process follows from
  \cite[Theorem 1]{Sag-91}, see \cite{KaS-98} for a proof. Specifically, we consider the limit along the subsequence $(\lambda_n,n\ge 1)$ for which 
  $\psi^{-1}(\lambda_n)=n$. 
\end{pf}
Proposition \ref{prop4a} states that the conclusion of Proposition \ref{Sagconv} holds in the supercritical as well as $\psi^\prime(\infty)<\infty$ cases. Sagitov \cite{Sag-91} announces his convergence result to include the (finite-mean) supercritical case, but the proof in Kaj and Sagitov \cite{KaS-98} only
treats the subcritical and critical cases. They say that the supercritical case would require some further assumptions and additional work. In our less general setting, this is not difficult -- our proof does not rely on \cite{KaS-98}.  

\begin{pfofprop4a} In the g.i.d.\ case $\psi^\prime(\infty)=\infty$, we use Lemma \ref{sub and geo} to represent $\kappa_\lambda$ in terms of a strictly increasing subordinator $\sigma$ as $\kappa_\lambda=\bP(\sigma(V_\lambda)\in\cdot)$ for 
  $V_\lambda\sim{\rm Exp}(\psi^\prime(\psi^{-1}(\lambda)))$. We introduce Markovian measure-valued branching processes $M^\lambda$ that do not record residual
  lifetimes of $Z^\lambda$, but what we may think of as limiting residual lifetimes (in the limit $\lambda\rightarrow\infty$),
  $$M^\lambda_t=\sum_{u\in\overline{\pi}_t(B_\lambda^\bullet)}\delta_{\sigma_u(K^u_{t-\alpha_u^\bullet})-(t-\alpha_u^\bullet)}=\sum_{x\in\cB^\circ_\lambda\colon\sigma(x-)\le t\le\sigma(x)}\delta_{\sigma(x)-t},\qquad\mbox{where $\sigma(x-)=\cW_\lambda([[0,x[[)$,}$$
%  although this is interpretation actually fails in the case $\psi^\prime(\infty)=\infty$ without affecting the validity of the argument. 
using notation as in Section \ref{secsnake} and also $K^u_s=\inf\{t\ge 0\colon\sigma_u(t)>s\}$.
By \cite[Lemma 3]{KaS-98} or \cite[Formula (1.5)]{Dyn-91}, for the Markov process $\xi_t=\sigma(K_t)-t$, the
  semigroup of $M^\lambda$ is such that
  $$\bE\left(\!\left.\exp\left\{\!-\!\int_{[0,\infty)}\!f(z)M_t^\lambda(dz)\right\}\right|M_0^\lambda\!=\!m\right)=\exp\left\{\!-\!\int_{[0,t]}\!v^\lambda_{t-z}(f)m(dz)+\int_{(t,\infty)}\!f(z\!-\!t)m(dz)\right\},$$
  for all $f\colon[0,\infty)\rightarrow[0,\infty)$ bounded continuous, where $v_t^\lambda(f)$ satisfies
  $$1-e^{-v_t^\lambda(f)}=\bE\left(1-e^{-f(\sigma(K_t)-t)}\right)-\int_{[0,t]}\left(\varphi_\lambda(e^{-v_{t-s}^\lambda(f)})-e^{-v_{t-s}^\lambda(f)}\right)d\bE(N_s^\lambda),$$
  and where $N_s^\lambda=\#\{k\ge 1\colon R_k^\lambda\le s\}$ is the renewal process associated with a random walk with $R_k^\lambda-R_{k-1}^\lambda\sim\kappa_\lambda$. It
  is easily checked that $\bE(N_s^\lambda)=H_s$. The remainder is straightforward (cf. \cite[Section 1.2]{Dyn-91}). We apply (\ref{psipsitilde}) to see that $u_t^\lambda(f)=\psi^{-1}(\lambda)(1-\exp\{-v_t^\lambda(f/\psi^{-1}(\lambda))\})$ satisfies
  $$u_t^\lambda(f)+\int_{[0,t]}\psi(u_{t-s}^\lambda(f))dH_s=\bE\left(\psi^{-1}(\lambda)\left(1-\exp\left\{-\frac{f(\sigma(K_t)-t)}{\psi^{-1}(\lambda)}\right\}\right)\right).$$
  Uniqueness in (\ref{measint}) means that $u_t^\lambda(f)=u_t(f_\lambda)$, where $f_\lambda=\psi^{-1}(\lambda)(1-e^{-f/\psi^{-1}(\lambda)})\uparrow f$. By 
  the Monotone Convergence Theorem, this implies for $N_\lambda\sim{\rm Poi}(\beta\psi^{-1}(\lambda))$ and $M_0^\lambda=N_\lambda\delta_0$ that
  \beq\bE\left(\!\exp\left\{\!-\!\!\int_{[0,\infty)}\!\!f(z)\frac{M_t^\lambda(dz)}{\psi^{-1}(\lambda)}\right\}\right)&\!\!\!=\!\!\!&\exp\left\{-\beta\psi^{-1}(\lambda)(1-v_t^\lambda(f/\psi^{-1}(\lambda)))\right\}=\exp\{-\beta u_t(f_\lambda)\}\\
       &\!\!\!=\!\!\!&\bE\left(\left.\!\exp\left\{\!-\!\!\int_{[0,\infty)}\!\!f_\lambda(z)M_t(dz)\!\right\}\right|M_0=\beta\delta_0\right)\\
       &\!\!\!\rightarrow\!\!\!&\bE\left(\left.\!\exp\left\{\!-\!\!\int_{[0,\infty)}\!\!f(z)M_t(dz)\!\right\}\right|M_0=\beta\delta_0\right)=\exp\{-\beta u_t(f)\}.
  \eeq
  In particular, for $Z_t^\lambda=M_t^\lambda([0,\infty)$, we obtain $Z_t^\lambda/\psi^{-1}(\lambda)\rightarrow Z_t$, where $Z$ is a ${\rm CSBP}(K,\psi)$. The martingale argument of Lemma \ref{CSBP Limit Lemma x} establishes almost sure convergence in the case $\psi^\prime(0)>-\infty$. 
  
  In the finitely geometrically divisible case, we use bushes $\overline{B}_\lambda$ based on measures $\bQ^{\rm RW}_\lambda$ and
  $$M_t^\lambda=\sum_{u\in\overline{\pi}({B}_\lambda^\bullet)}\delta_{\sigma_u(K^u_{t-\alpha_u^\bullet})-(t-\alpha_u^\bullet)},\qquad\mbox{where $K^u_s=\inf\{k\ge 1\colon \sigma_u(k)>s\}$.}$$
  Then the argument above is easily adapted.
\end{pfofprop4a}

As an application, let us derive a backbone decomposition. This should be useful to deduce more general supercritical Bellman-Harris convergence results from subcritical results. Our present paper is not about superprocesses nor convergence of general triangular arrays, so we do not push for highest generality nor assumptions as in \cite{KaS-98}, but we would like to mention a now natural approach -- in a sense to be made precise, convergence of supercritical processes is equivalent to convergence of backbones and convergence of associated subcritical processes.

We write $\bP_{K,\psi}^r$ for the distribution of a ${\rm CSBP}(K,\psi)$ starting from $r\ge 0$. Just as for ${\rm CSBP}(\psi)$ in Section \ref{CSBP and CBI explanation
section}, we consider the sigma-finite measure $\Theta_{K,\psi}$ such that $\bP_{K,\psi}^r$ is the distribution of a sum over a Poisson point process with intensity measure $r\Theta_{K,\psi}$. 

\begin{theo}[Backbone decomposition for ${\rm CSBP}(K,\psi)$]\label{CSBPbackbone} Let $\psi$ be a (non-explosive) supercritical, $\psi_0(r)=\psi(r+\psi^{-1}(0))$ the associated subcritical branching mechanism. Let $\overline{B}_0$ be a bush of 
  $N_0\sim{\rm Poi}(\beta\psi^{-1}(0))$ trees with distribution 
  $\bP_{\bQ_0}$ as in {\rm (\ref{Qmu})}, and, as in Section \ref{secsnake}, $(\cB_0^\circ,\cW_0)$ a representation as a weighted $\bR$-tree, $\sigma(x)=\cW_0([[0,x]])$. Given $(\cB_0^\circ,\cW_0)$, consider
  \begin{itemize}\item points $(Z^x,x\in\cP)$ of a Poisson point process in $\bD([0,\infty),[0,\infty))$ with intensity measure
    $$Q(df,dx)=\left(2a\Theta_{K,\psi_0}(df)+\int_{(0,\infty)}\bP^r_{K,\psi_0}(df)re^{-r\psi^{-1}(0)}\Pi(dr)\right){\rm Leb}|_{\cB_0^\circ}(dx),$$
    \item extra points $(Z^x,x\in{\rm Br}(\cB_0^\circ))$ independent of $(Z^x,x\in\cP)$ with distribution
    $$Q^{(l(x))}(df)=\frac{2a1_{\{l(x)=2\}}}{|\psi_0^{(l(x))}(0)|}\delta_0(df)+\int_{(0,\infty)}\bP^r_{K,\psi_0}(df)\frac{r^{l(x)}e^{-r\psi^{-1}(0)}}{|\psi^{(l(x))}_0(0)|}\Pi(dr),$$
    where $l(x)+1$ is the number of connected components of $\cB_0^\circ\setminus\{x\}$ and ${\rm Br}(\cB_0^\circ)$ is the set of branchpoints
     $\{x\in\cB_0^\circ\colon l(x)\ge 2\}\setminus\{0\}$, 
    \item and an extra point $Z^{0}$ independent of $(Z^x,x\in\cP\cup{\rm Br}(\cB_0^\circ))$ with distribution $\bP^\beta_{K,\psi_0}$.
  \end{itemize}
  Then the process $\displaystyle Z_t=\sum_{x\in\cP\cup{\rm Br}(\cB_0^\circ)\cup\{0\}}Z^x_{t-\sigma(x)}$ is a ${\rm CSBP}(K,\psi)$ starting from $Z_0=\beta$.
\end{theo}
\begin{pf} Since $Z$ is not a Markov process, we will deduce the theorem from the richer structure of a Markovian $(\xi,K,\psi)$-superprocesses $M$ starting from
  $M_0=m$. Slightly abusing notation, we consider $M$ also under $\Theta_{K,\psi}$ and $\bP^r_{K,\psi}$ and $Q$. Then the intensity measures and distributions in the   
  bullet points specify a point process $(M^x,x\in\cP\cup{\rm Br}(\cB_0^\circ)\cup\{0\})$. From the exponential formula for
  Poisson point processes and from (\ref{meassg}) for the subcritical branching mechanism $\psi_0$ with $u^{\rm sub}_t(f)$ associated via the analogue of (\ref{measint}), it 
  is not hard to calculate 
  \begin{eqnarray*}&&\hspace{-0.68cm}\bE\left(\!\left.\exp\left\{\!-\!\sum_{x\in\cP}\int_{[0,\infty)}\!\!f(z)M^x_{t-\sigma(x)}(dz)\!\right\}\right|\cB_0^\circ,\cW_0\!\right)\!=\exp\left\{\!-\!\int_{\cB_0^\circ}\!(\psi_0^\prime(u_{t-\sigma(x)}^{\rm sub}(f))\!-\!\psi_0^\prime(0)){\rm Leb}(dx)\!\right\},\\
	  &&\hspace{-0.68cm}\bE\left(\left.\exp\left\{-\sum_{x\in{\rm Br}(\cB_0^\circ)}\int_{[0,\infty)}f(z)M^x_{t-\sigma(x)}(dz)\right\}\right|\cB_0^\circ,\cW_0\right)=\prod_{x\in{\rm Br}(\cB_0^\circ)}\frac{\psi_0^{(l)}(u_{t-\sigma(x)}^{\rm sub}(f))}{\psi_0^{(l)}(0)},\\
      &&\hspace{-0.68cm}\bE\left(\!\left.\exp\left\{\!-\!\int_{[0,\infty)}\!f(z)M^0_t(dz)\!\right\}\right|\cB_0^\circ,\cW_0\!\right)\!=\exp\left\{\!-\!\int_{[0,t]}\!u_{t-z}^{\rm sub}(f)m(dz)\!-\!\int_{(t,\infty)}\!f(z\!-\!t)m(dz)\!\right\},
  \end{eqnarray*}
  using the convention $u^{\rm sub}_s(f)=0$ for $s<0$. It now suffices to show that the backbone decomposition of approximations $M^\lambda$ of $M$, cf. Corollary 
  \ref{cor14}, appropriately converges to these quantities. Let us formulate that backbone decomposition in the current setting. Let  
  $\bP^{r,\lambda}_{\kappa,\psi_0}$ be the distribution of $M^\lambda$ given $M_0^\lambda=N_\lambda^r\delta_0$ with $N_\lambda^r\sim{\rm Poi}(\psi^{-1}(\lambda)r)$. Then, given $(\cB_0^\circ,\cW_0)$, consider points 
  $(M^{x,\lambda},x\in\cP^\lambda)$ of a Poisson point process with intensity measure 
  $$2a(\psi^{-1}(\lambda)-\psi^{-1}(0))\bP(M^\lambda\in\cdot|M^\lambda_0=\delta_0)+\int_{(0,\infty)}\bP^{r,\lambda}_{\kappa,\psi_0}re^{-r\psi^{-1}(0)}\Pi(dr),$$
  and  
$$M^{x,\lambda}\sim\frac{2a1_{\{l(x)=2\}}}{|\psi_0^{(l(x))}(0)|}\delta_0+\int_{(0,\infty)}\bP^{r,\lambda}_{K,\psi_0}\frac{r^{l(x)}e^{-r\psi^{-1}(0)}}{|\psi^{(l(x))}_0(0)|}\Pi(dr)\quad\mbox{for $x\in{\rm Br}(\cB_0^\circ)$ and } M^{0,\lambda}\sim\bP^{\beta,\lambda}_{K,\psi_0}.$$
  Then the analogous calculations yield for $u^{\lambda,\rm sub}_t(f)=\psi_0^{-1}(\lambda)(1-e^{-v_t^{\lambda,\rm sub}(f/\psi^{-1}(\lambda))})$
 \begin{eqnarray*}&&\hspace{-0.68cm}\bE\left(\!\left.\exp\left\{\!-\!\!\sum_{x\in\cP^\lambda}\int_{[0,\infty)}\!\!\!\!\!f(z)\frac{M^{x,\lambda}_{t-\sigma(x)}(dz)}{\psi^{-1}(\lambda)}\!\right\}\right|\cB_0^\circ,\cW_0\!\right)\!=\exp\left\{\!-\!\int_{\cB_0^\circ}\!\!(\psi_0^\prime(u^{\lambda,\rm sub}_{t-\sigma(x)}(f))\!-\!\psi_0^\prime(0)){\rm Leb}(dx)\!\right\}\!,\\
	  &&\hspace{-0.68cm}\bE\left(\left.\exp\left\{-\sum_{x\in{\rm Br}(\cB_0^\circ)}\int_{[0,\infty)}f(z)\frac{M^x_{t-\sigma(x)}(dz)}{\psi^{-1}(\lambda)}\right\}\right|\cB_0^\circ,\cW_0\right)=\prod_{x\in{\rm Br}(\cB_0^\circ)}\frac{\psi_0^{(l)}(u_{t-\sigma(x)}^{\lambda,\rm sub}(f))}{\psi_0^{(l)}(0)},\\
      &&\hspace{-0.68cm}\bE\left(\!\left.\exp\left\{\!-\!\int_{[0,\infty)}\!f(z)\frac{M^0_t(dz)}{\psi^{-1}(\lambda)}\!\right\}\right|\cB_0^\circ,\cW_0\!\right)\!=\exp\left\{\!-\!\int_{[0,t]}\!u_{t-z}^{\lambda,\rm sub}(f)m(dz)\!-\!\int_{(t,\infty)}\!f(z\!-\!t)m(dz)\!\right\},
  \end{eqnarray*}
  But from the proof of Proposition \ref{prop4a}, we know that $u_t^{\lambda,\rm sub}(f)\rightarrow u_t^{\rm sub}(f)$ as $\lambda\rightarrow\infty$, and this completes the proof.  
\end{pf} 

We can specialise this backbone decomposition to the case $K_t=t$, when a ${\rm CSBP}(K,\psi)$ is simply a Markovian ${\rm CSBP}(\psi)$ and $\sigma(x)=d(0,x)$ is just the height of $x\in\cB_0^\circ$. In this framework, and even with a spatial motion added, this
decomposition was obtained recently by Berestycki et al. \cite{BKM-10}, generalising an analogous result of \cite[Theorem 5.6]{DuW} formulated in a context of L\'evy trees.

%-----------------------------------------------------------------------------------------------------------------------------------Limit to CSBP

%-----------------------------------------------------------------------------------------------------------------------------------------Immigrations
\section{Growth of ${\rm GWI}(q_\lambda,\kappa_\lambda,\eta_\lambda,\chi_\lambda)$-forests: immigration}\label{Section Immigration}

Theorem \ref{thm3} is about forests $F_\lambda$ of ${\rm GW}(q_\lambda,{\rm Exp}(c_\lambda))$-trees arising from immigration of independent $\eta_\lambda$-distributed
numbers of immigrants at ${\rm Exp}(\paramgamma_\lambda)$-spaced times. The main statement beyond the no-immigration case of Theorem \ref{thm1} is that consistency of a
family $(F_\lambda,\lambda\ge 0)$ under Bernoulli leaf colouring relates $\eta_\lambda$ to a continuous-state immigration
mechanism $\phi$. After some more general remarks, we focus on the Markovian case of exponential lifetimes and inter-immigration times.

\subsection{A two-colours regenerative property and associated forest reduction}\label{regprop}

Let $F=(B(t),t\ge 0)$ be a ${\rm GWI}(q,\kappa,\eta,\chi)$-forest as defined in Section \ref{discr BP with immigration section}, specifically denote by $S_1=\inf\{t\ge 0\colon B(t)\neq\partial\}\sim\chi$ the first immigration time, by $B(S_1)=(T_{(1)}^{(1)},\ldots,T_{(N_1)}^{(1)})$ the 
bush of independent genealogical trees $T^{(1)}_{(j)}\sim{\rm GW}(q,\kappa)$, $j\ge 1$, of the $N_1\sim\eta$ time-$S_1$ immigrants and by $F_{\rm post}=(B^{\rm post}(t),t\ge 0)$ the post-$S_1$ forest
given by $B^{\rm post}(0)=\partial$ and $B^{\rm post}(t)=B(S_1+t)$ for $t>0$. It is immediate from the definition that $F$ satisfies a
\em regenerative property \em at $S_1$ in that $(S_1,B(S_1))$ is independent of $F_{\rm post}$ and $F_{\rm post}\overset{{\rm (d)}}{=}  F$, and that the distribution of $(S_1,B(S_1))$ as above together with this regenerative
property characterises the distribution of $F$. Since colouring and reduction apply tree by tree, we obtain for the associated
forest $F^{p-\rm col}$ of coloured trees $T^{(i),p-\rm col}_{(j)}\sim\bP^{p-\rm col}_{q\otimes\kappa}$, 
$1\le j\le N_i$, $i\ge 1$, a

\paragraph{Regenerative property of coloured ${\rm GWI}(q,\kappa,\eta,\chi)$-forests}
\begin{enumerate}\item[(a)] For all $n\ge 1$, $\varepsilon_j\in\{0,1\}$, and measurable functions $k$, $f_j$ and $G$, $1\le j\le n$,
  we have
  \begin{eqnarray*}&&\hspace{-1.5cm}\bE\left(k(S_1)G(F^{p-\rm col}_{\rm post})\prod_{j=1}^nf_j(T^{p-\rm col}_{(j)});N_1=n;\left(\gamma_\emptyset(T^{p-\rm col}_{(1)}),\ldots,\gamma_\emptyset(T^{p-\rm col}_{(n)})\right)=(\varepsilon_1,\ldots,\varepsilon_n)\right)\\
   &&=\int_{(0,\infty)}k(z)\chi(dz)\eta(n)g(p)^{n_r}(1-g(p))^{n_b}\bE(G(F^{p-\rm col}))\prod_{j=1}^n\bE_{q\otimes\kappa}^{p-\rm col}[f_j|\gamma_\emptyset=\varepsilon_j],
  \end{eqnarray*}
  where $n_r=\varepsilon_1+\cdots+\varepsilon_n$ and $n_b=n-n_r$ are the numbers of red and black colours.
  \item[(b)] For $t\ge 0$, consider the post-$t$ forest $F_{{\rm post}-t}^{p-\rm col}=(B^{p-\rm col}(t+s),s\ge 0)$ and 
    the pre-$t$ sigma-algebra $\cF_t=\sigma\{B^{p-\rm col}(r),r\le t\}$. Then for all measurable functions $f_u$, $u\in\bU$, and \nolinebreak $G$,
  $$\hspace{-0.5cm}\bE\!\left(\!\left.G(F_{{\rm post}-t}^{p-\rm col})\!\!\!\prod_{u\in\overline{\pi}_t(F)}\!\!f_u(\overline{\theta}_{u,t}(F^{p-\rm col}))\right|\cF_t\right)\!
	=\bE(G(F^{p-\rm col}))\!\!\!\prod_{u\in\overline{\pi}_t(F)}\!\!\left.\bE^{p-\rm col}_{q\otimes\kappa}[f_u(\overline{\theta}_{\emptyset,s})|\zeta_\emptyset>s]\right|_{s=t-\alpha_u}\!.$$
\end{enumerate}
As a trivial application of (a), we can calculate the probability that all immigrants are red
$$\bP\left(\left(\gamma_\emptyset(T^{p-\rm col}_{(1)}),\ldots,\gamma_\emptyset(T^{p-\rm col}_{(N_1)})\right)=(1,\ldots,1)\right)
=\sum_{n\ge 1}\eta(n)g(p)^n=\varphi_{\eta}(g(p)).
$$
We deduce the distribution of the number of red immigrants given that all immigrants are red
\begin{equation}\label{red disbribution}
 \eta_{\rm red}^{p-\rm col}(m)=\left\{
  \begin{array}{ll}
    \frac{\eta(m)g(p)^{m}}{\varphi_{\eta}(g(p))}, & \hbox{if $m \geq 1$,} \\
    0, & \hbox{if $m=0$,}
  \end{array}
\right.\quad\mbox{with generating function }\varphi_{\eta_{\rm red}^{p-\rm col}}(s)=\frac{\varphi_\eta(sg(p))}{\varphi_\eta(g(p))}.
\end{equation}
Also by the regenerative property (a), the number $\widetilde{G}$ of immigrations until we see the first black immigrant is geometrically distributed with parameter $1-\varphi_\eta(g(p))$, i.e.
\begin{equation}\label{Immigration Geometric}
\bP(\widetilde{G}=j)=\varphi_{\eta}(g(p))^{j-1}(1-\varphi_{\eta}(g(p))),\qquad j\ge 1.
\end{equation}
Conditioning on having at least one black immigrant (probability $1-\varphi_\eta(g(p))$), we get for $\ell\ge 1$
$$\eta^{p-\rm rdc}(\ell)=\frac{1}{1-\varphi_\eta}(g(p))\sum_{m\ge 0}{m+\ell\choose m}\eta(m+\ell)(g(p))^m(1-g(p))^\ell=\frac{(1-g(p))^\ell\varphi_\eta^{(\ell)}(g(p))}{\ell!(1-\varphi_\eta(g(p))}.$$
Similarly, conditioning on having $\ell$ black immigrants, $\ell\ge 1$, we get for the number of red ones
\begin{equation}\label{red conditional distribution}\widetilde{\eta}_\ell(m)=\eta(m+\ell)\frac{(m+\ell)!}{m!}(g(p))^m\frac{1}{\varphi_\eta^{(\ell)}(g(p))},\qquad\mbox{for $m\ge 0$.}\end{equation}
These distributions have generating functions that we can express in terms of $\varphi_\eta$, for $s\in[0,1]$
$$\varphi_{\eta^{p-\rm rdc}}(s)=\frac{\varphi_\eta(g(p)+s(1-g(p)))-\varphi_\eta(g(p))}{1-\varphi_\eta(g(p))}\qquad\mbox{and}\qquad\varphi_{\widetilde{\eta}_\ell}(s)=\frac{\varphi_\eta^{(\ell)}(sg(p))}{\varphi_\eta^{(\ell)}(g(p))},\quad\ell\ge 1,$$
and as $\varphi_{\eta^{p-\rm rdc}}$ and $\varphi_\eta$ are analytic, we can extend $\varphi_\eta$ analytically to $[-g(p)/(1-g(p)),1]$. Evaluating at $s=v_p=-g(p)/(1-g(p))$, we get
$\varphi_\eta(g(p))=-\varphi_{\eta^{p-\rm rdc}}(v_p)/(1-\varphi_{\eta^{p-\rm rdc}}(v_p))$, so
\begin{equation}\label{exteta}\varphi_\eta(r)=\frac{\varphi_{\eta^{p-\rm rdc}}(v_p+r(1-v_p))-\varphi_{\eta^{p-\rm rdc}}(v_p)}{1-\varphi_{\eta^{p-\rm rdc}}(v_p)},\qquad r\in[0,1].
\end{equation} 
As reduction preserves the Galton-Watson property for trees \cite{DuW}, we now see that a $p$-reduced ${\rm GWI}(q,\kappa,\eta,\chi)$-forest is a ${\rm GWI}(q^{p-\rm rdc},\kappa^{p-\rm rdc},\eta^{p-\rm{rdc}}, \chi^{p-\rm{rdc}})$-forest. Specifically, $\eta^{p-\rm rdc}$ is as above, $\chi^{p-\rm{rdc}}$ the distribution of a geom($1-\varphi_\eta(g(p))$) sum of independent
$\chi$-distributed variables.

\subsection{Growth of ${\rm GWI}(q_\lambda,{\rm Exp}(c_\lambda),\eta_\lambda,{\rm Exp}(\paramgamma_\lambda))$-forests}

\subsubsection{Proof of Theorem \ref{thm3}}
%----------------------------reasoning for immigrations-------------------------------------------------------------------------------------
${\rm (i)}\Rightarrow{\rm (ii)}$: Suppose, (i) holds. In particular, $q_\mu$ is then the $(1-\mu/\lambda)$-reduced
offspring distribution associated with $q_\lambda$, for all $0\le\mu<\lambda<\infty$. By Theorem \ref{thm1},
$\varphi_q(s)=s+\widetilde{\psi}(1-s)$, where $\widetilde{\psi}$ has the form (\ref{lkbm}). To be specific, let $c=1$ and parametrise $(q_\lambda,c_\lambda)$ using $\psi$ as in Section \ref{Consistent family}.

Now consider the relationship between $\eta=\eta_\mu$ and $\eta_\lambda$ for $\lambda>\mu=1$. By the discussion above (\ref{exteta}), we can extend $\varphi_\eta$ analytically to $[-g_\lambda(1-1/\lambda)/(1-g_\lambda(1-1/\lambda)),1]$, where, expressing as in (\ref{g in psi}), we have $g_\lambda(1-1/\lambda)=1-\psi^{-1}(1)/\psi^{-1}(\lambda)\rightarrow 1$ as $\lambda\rightarrow\infty$. Differentiating (\ref{exteta}),
we see that $\varphi_\eta^\prime$ has positive derivatives on $(-\infty,1)$. Setting $\widetilde{\phi}(r)=1-\varphi_\eta(1-r)$, 
$r\ge 0$, the derivative $\widetilde{\phi}^\prime$ is completely monotone on $(0,\infty)$ and, by Bernstein's theorem (see e.g. \cite{Fel}), there exists a
Radon measure $\widetilde{\Lambda}^*$ on $[0, \infty)$ such that
$$\widetilde{\phi}^{\prime}(r)=\int_{[0,\infty)}e^{-r x} \widetilde{\Lambda}^*(d x)<\infty,\qquad r>0.$$
From $\widetilde{\phi}(0)=0$, we get integrability $\int_{(1,\infty)}x^{-1}\widetilde{\Lambda}^*(dx)<\infty$ and 
$$
  \widetilde{\phi}(u) = \widetilde{\phi}(0) + \int^{u}_{0} \widetilde{\phi}^{\prime}(r) d r
                = \widetilde{\Lambda}^*(\{0\})u + \int_{(0,\infty)}\frac{1-e^{-ux}}{x}\widetilde{\Lambda}^*(d x),
$$
and, in particular, setting $\widetilde{d}=\widetilde{\Lambda}^*(\{0\})$ and $\widetilde{\Lambda}(dx)=x^{-1}\widetilde{\Lambda}^*|_{(0,\infty)}(dx)$ yields (ii).

%-----------------------------------------------------------------------------------------------------------------------------------------------
${\rm (ii)}\Rightarrow{\rm (i)}$: Now suppose that (ii) holds. According to Theorem \ref{thm1}, the family of offspring distribution $(q_\lambda,\lambda\ge 0)$ exists as required; furthermore, we can express $\varphi_{q_\lambda}$ in terms of $\psi$ as in Section \ref{Consistent family}, choosing $c=1$.

By (\ref{g in psi}), we have $g_\lambda(1-\mu/\lambda)=1-\psi^{-1}(\mu)/\psi^{-1}(\lambda)$, so $v_{\mu,\lambda}=-g_\lambda(1-\mu/\lambda)/(1-g_\lambda(1-\mu/\lambda))=1-\psi^{-1}(\lambda)/\psi^{-1}(\mu)$ for all $0\le\mu<\lambda<\infty$.
By (\ref{exteta}) and the discussion above (\ref{exteta}), the required immigration distributions must be of the following form, respectively for $\mu<1$ and $\lambda>1$
\beq\varphi_{\eta_\mu}(s)=\frac{\varphi_\eta(g_1(1-\mu)+s(1-g_1(1-\mu)))-\varphi_\eta(g_1(1-\mu))}{1-\varphi_\eta(g_1(1-\mu))}\!\!&=&\!\!1-\frac{\widetilde{\phi}((1-s)\psi^{-1}(\mu)/\psi^{-1}(1))}{\widetilde{\phi}(\psi^{-1}(\mu)/\psi^{-1}(1))},\\
\varphi_{\eta_\lambda}(r)=\frac{\varphi_\eta(v_{1,\lambda}+r(1-v_{1,\lambda}))-\varphi_\eta(v_{1,\lambda})}{1-\varphi_\eta(v_{1,\lambda})}\!\!&=&\!\!1-\frac{\widetilde{\phi}((1-r)\psi^{-1}(\lambda)/\psi^{-1}(1)}{\widetilde{\phi}(\psi^{-1}(\lambda)/\psi^{-1}(1))}.\eeq
Since $\widetilde{\phi}^\prime$ is completely monotone, simple differentiation yields that these functions are indeed generating functions of immigration distributions. Furthermore, $\eta_\mu$ is the $(1-\mu/\lambda)$-reduced immigration 
distribution of $\eta_\lambda$ for all $0\le\mu<\lambda<\infty$, by the transitivity of colouring reduction noted in
Remark \ref{remdet}(a), which also applies to forests, since colouring and reduction are defined tree by tree. The full 
statement of (i) can now be obtained formally as in the proof of Theorem \ref{thm2}, with the simpler regenerative property here taking the role of the branching property there.

In the setting of (i) and (ii) for $c=c_1\in(0,\infty)$, $\paramgamma=\paramgamma_1\in(0,\infty)$, Kolmogorov's consistency theorem allows us to set up a consistent family $(F_\lambda,\lambda\ge 0)$ of ${\rm GWI}(q_\lambda,{\rm Exp}(c_\lambda),\eta_\lambda,{\rm Exp}(\paramgamma_\lambda))$-forests. Uniqueness of ($q_\lambda,c_\lambda$), $\lambda \geq 0$,
follows from Theorem \ref{thm1}. Uniqueness of ($\eta_\lambda,\lambda\ge 0$) was noted in ${\rm (ii)}\Rightarrow{\rm (i)}$.
Uniqueness of $(\paramgamma_\lambda,\lambda\ge 0)$ follows from the relationship between inter-immigration times as geometric sums,
where we calculate $\paramgamma_\lambda$ from (\ref{Immigration Geometric}) for $\lambda>1>\mu$ as 
$$\paramgamma_\lambda=(1-\varphi_{\eta_\lambda}(g_\lambda(1-1/\lambda)))\paramgamma=\frac{\paramgamma}{\widetilde{\phi}(\psi^{-1}(\lambda)/\psi^{-1}(1))}\quad\mbox{and}\quad \paramgamma_\mu
%=\frac{\paramgamma}{1-\varphi_{\eta_1}(g_1(1-\mu))}
=\frac{\paramgamma}{\widetilde{\phi}(\psi^{-1}(\mu)/\psi^{-1}(1))}.$$
%\vspace{-0.6cm}
\begin{flushright}$\Box$\end{flushright}

\subsubsection{Freedom in parameterisation and standard choice}

In analogy to Section \ref{Consistent family}, we can use a single function $\phi$ to replace $(\widetilde{\phi},\paramgamma)$ of Theorem \ref{thm3} and parametrise $(\eta_\lambda,{\rm Exp}(\paramgamma_\lambda))$, $\lambda\ge 0$, such that
\begin{equation}\label{full immi gen by phi psi}
\varphi_{\eta_\lambda}(v)=1-\frac{{\phi}(\psi^{-1}(\lambda)(1-v))}{{\phi}(\psi^{-1}(\lambda))}\qquad\mbox{and}\qquad\paramgamma_\lambda=\phi(\psi^{-1}(\lambda)),
\end{equation}
where $\phi$ is a linear transformation $\phi(s)=k_3\widetilde{\phi}(k_4s)$ of $\widetilde{\phi}$. Specifically, we choose $k_3=\paramgamma$ and $k_4=1/\psi^{-1}(1)$. It is easy to check that this works, using $\widetilde{\phi}(1)=1-\eta(0)=1$. 

In this parameterisation, we can also express in terms of $\psi$ and $\phi$ the remaining quantities studied in Section \ref{regprop}. E.g. (\ref{full immi gen by phi psi}) and (\ref{red disbribution}) now yield the generating function of the pure-red 
immigration distribution 
\begin{equation}\label{red immi gen}
    \varphi_{\eta_{\lambda,\rm red}^{(1-\mu/\lambda)-\rm col}}(s)=\frac{\phi(\psi^{-1}(\lambda))-\phi( \psi^{-1}(\lambda)-s(\psi^{-1}(\lambda)-\psi^{-1}(\mu)))}{\psi^{-1}(\lambda)-\psi^{-1}(\mu)}.
\end{equation}
The parameters of the geometric distributions (\ref{Immigration Geometric}) take a simple form that leads to a distinction of finite/infinite immigration rate $\alpha^{\rm imm}_{\mu,\lambda}=1-\varphi_{\eta_\lambda}(g_\lambda(1-\mu/\lambda))=\phi(\psi^{-1}(\mu))/\phi(\psi^{-1}(\lambda))$, and 
\begin{itemize}
\item $\alpha_{\mu,\lambda}^{\rm imm} \rightarrow\frac{\phi(\psi^{-1}(\mu))}{\phi(\infty)}>0$ as $\lambda \rightarrow \infty$, if $\phi(\infty) < \infty;$
\item $\alpha_{\mu,\lambda}^{\rm imm} \rightarrow 0$ as $\lambda \rightarrow \infty$, if $\phi(\infty) = \infty$.
\end{itemize}
With the formulas above we can formulate explicitly a reconstruction procedure.

\subsubsection{Reconstruction procedure for ${\rm GWI}(q_\lambda,{\rm Exp}(c_\lambda),\eta_\lambda,{\rm Exp}(\paramgamma_\lambda))$-forests}

For $F_\mu\sim{\rm GWI}(q_\mu,{\rm Exp}(c_\mu),\eta_\mu,{\rm Exp}(\paramgamma_\mu))$, we modify the steps of Section 
\ref{Recon for GW expo bushes2} to construct $F_\lambda$.
\begin{enumerate}\item[1.] In every tree of $F_\mu$, subdivide lifetimes as in Section \ref{Recon for GW expo bushes2} and hence
  construct a forest $\widehat{F}_\mu$.
  \item[2.] In every tree of $\widehat{F}_\mu$, add further children and independent red trees as in Section
    \ref{Recon for GW expo bushes2} and hence construct a forest $\widehat{F}_\lambda$.
  \item[3.] At every immigration time, given that there are $N_\mu^{(i)}=\ell$ immigrants in $\widehat{F}_\mu$, consider a
    random number $N_\lambda^{(i),\rm red}\sim\widetilde{\eta}_\ell$ of further immigrants as in (\ref{red conditional distribution}),
    proceed as in Section \ref{Recon for GW expo bushes2} and superpose a further independent ${\rm GWI}(q_{\lambda,{\rm red}}^{(1-\mu/\lambda)-\rm col},{\rm Exp}(c_\lambda),\eta_{\lambda,\rm red}^{(1-\mu/\lambda)-\rm col},{\rm Exp}(\paramgamma_\lambda-\paramgamma_\mu))$-forest with distributions as in (\ref{red offspring distribution}) and (\ref{red immi gen}) to finally obtain $F_\lambda$.
\end{enumerate}

%------------------------------------------------------------------------------------------------------------------------------------------
%---------------------
%------------------------------
%-----------------------------------------------------------------------------------------------------------------------Poisson Immigration Reconstruction.

%-------------------------------------------------------------------------------------------------Immigration and Edges!!!!!!!!
\subsubsection{Convergence of the population sizes: proof of Proposition \ref{prop4}}\label{section convergence of the population size}

For convergence in distribution, we calculate the Laplace transform of $Y_t^\lambda/\psi^{-1}(\lambda)$, where
$$Y_t^\lambda=\overline{\pi}_t\left(F_\lambda\right)=\sum_{i=1}^{J_t^\lambda}\overline{\pi}_{t-S_i^\lambda}\left(B_\lambda(S_i^\lambda)\right)=\sum_{i=1}^{J_t^\lambda}\sum_{j=1}^{N_i^\lambda}\overline{\pi}_{t-S_i^\lambda}\left(T_{(j)}^{(i),\lambda}\right)$$
with notation as in and around (\ref{discrimm}), but with all quantities $\lambda$-dependent.
We exploit that
$$E^\lambda(S_i^\lambda)=Z^{(i),\lambda}=\left(\overline{\pi}_t(B_\lambda(S_i^\lambda)),t\ge 0\right),\qquad E^\lambda(s)=0,\quad s\not\in\left\{S_i^\lambda,i\ge 1\right\},$$
is a Poisson point process with intensity measure $\paramgamma_\lambda$ times
the distribution of a ${\rm GW}(q_\lambda,{\rm Exp}(c_\lambda))$-process $Z^{(1),\lambda}$ starting from $Z^{(1),\lambda}_0\sim\eta_\lambda$. By the exponential formula for Poisson point processes,
$$\bE\left(\exp\left\{-rY_t^\lambda/\psi^{-1}(\lambda)\right\}\right)=\exp\left\{-\paramgamma_\lambda\int_0^t\sum_{m=1}^\infty\eta_\lambda(m)\left(1-\left(\bE\left(s^{Z_{t-v}^\lambda}\right)\right)^m\right)dv\right\},$$
where $s=e^{-r/\psi^{-1}(\lambda)}$ and $Z^\lambda$ is the population size of a single ${\rm GW}(q_\lambda,{\rm Exp}(c_\lambda))$-tree as in the proof of Lemma \ref{CSBP Limit Lemma x}. Using notation and asymptotics from there, as well as (\ref{full immi gen by phi psi}), this equals
$$\exp\left\{-\paramgamma_\lambda\int_0^t\left(1-\varphi_{\eta_\lambda}\left(w_{t-v}^\lambda(s)\right)\right)dv\right\}=\exp\left\{-\int_0^t\phi\left(\psi^{-1}(\lambda)\left(1-w_{t-v}^\lambda(s)\right)\right)dv\right\}.$$
Since $\psi^{-1}(\lambda)\left(1-w_{t-v}^\lambda(s)\right)\rightarrow u_{t-v}(r)$, and, by (\ref{CSBP u bounds}) and (\ref{wlambda}), all these quantities are bounded by
$\max\{r,\psi^{-1}(0)\}$, dominated convergence completes the proof of convergence in distribution.

Almost sure convergence follows by martingale arguments as in Lemma \ref{CSBP Limit Lemma x}, using a version of the regenerative property (ii) of Section \ref{regprop} rather than the version of the branching property (\ref{bpexp}) that we presented in Remark \ref{Markov}(c). {\hspace*{\fill} $\square$}
% can also show the almost surely convergence here. Recall the notation: 
%%$\overline{\pi}_t(T)\subset\bU$ for the population of tree $T$ alive at time $t$, and 
%the sigma-algebra
%$\mathcal{G}_{\lambda, t}=\sigma\left(\overline{\pi}_s(B_{\lambda^\prime}),\lambda^\prime \geq \lambda, s \leq t \right)$ 
%%and $\widetilde{\mathcal{G}}_{\lambda, t}=\sigma\left(\overline{\pi}_s(B_\lambda), s \leq t \right).$
%It can be shown that 
%$$\mathbb{E}\left(\frac{1}{\psi^{-1}(\mu)}Y^{\mu}_t | \mathcal{G}_{\lambda, t}\right)= %\frac{1}{\psi^{-1}(\lambda)}Y^{\lambda}_t,$$
%which is the martingale property in the (decreasing) filtration $(\mathcal{G}_{\lambda,t},\lambda\ge 0)$ that implies %$Y^{\lambda}_{t}/\psi^{-1}(\lambda) \rightarrow Y_t$ almost surely as $\lambda \rightarrow
%\infty$.

\bigskip

From Proposition \ref{prop4} and Lemma \ref{CSBP Limit Lemma x} we deduce the analogous convergence result for ${\rm GWI}$-processes starting from initial
population sizes $Y_0^\lambda\sim{\rm Poi}(\beta_\lambda)$. The limiting CBI then has $Y_0=\beta$. 
%
%\begin{coro}\label{Imm Limit Thm general}
%\rm Let $(Y^{\lambda, x}_t)_{t>0}$ denote the population size
%process in the setting of Theorem \ref{thm3} and $Y^{\lambda,x}_0=x \geq 0$, for any time $t
%\geq 0$ and $\psi'(0) > - \infty,$, a.s.
%$$\frac{1}{\psi^{-1}(\lambda)}Y^{\lambda,x}_t \stackrel{\lambda \rightarrow \infty}{\longrightarrow} Y^x_t$$
%\end{coro}
%\textbf{Proof:} Set $Y^{\lambda, x}_0$ into two parts
%$Z^{\lambda,x}_0$ and $Y^{\lambda, 0}_0.$ Firstly we can show that
%as $\lambda \rightarrow \infty,$
%$$\frac{Y^{\lambda, x}_0}{\psi^{-1}(\lambda)} \sim \frac{Z^{\lambda,x}_0}{\psi^{-1}(\lambda)} \rightarrow Z^{x}_0 \sim %Y^{x}_0,$$
%Then for $t > 0,$
%\begin{eqnarray*} \mathbf{E}(e^{-\mu
%\frac{Y^{\lambda, x}_t}{\psi^{-1}(\lambda)}}) &\rightarrow& \exp\{-x
%u_t(\mu)-\int^{t}_{0} \phi (u_r(\mu))ds\}\\
%&=& \mathbf{E}(e^{-\mu Y^{x}_t})
%\end{eqnarray*}
%\begin{flushright}$\Box$\end{flushright}

\subsection{Analogous results for ${\rm GWI}(q_\lambda,\kappa_\lambda,\eta_\lambda,\chi_\lambda)$-forests}

Finally, let us combine Theorems \ref{thm2} and \ref{thm3} into a single statement and also deduce the analogous pattern for general inter-immigration distributions that now emerges naturally.

\begin{coro}\label{corolast} For a tuple $(q,\kappa,\eta,\chi)$ of offspring, lifetime, immigration and inter-immigration distributions, the
  following are equivalent:
  \begin{enumerate}\item[\rm(i)] There are $(q_\lambda,\kappa_\lambda,\eta_\lambda,\chi_\lambda)_{\lambda\ge 0}$
    with $(q_1,\kappa_1,\eta_1,\chi_1)=(q,\kappa,\eta,\chi)$ such that
    $(q_\mu,\kappa_\mu,\eta_\mu,\chi_\mu)$ is the $(1-\mu/\lambda)$-reduced tuple associated with
    $(q_\lambda,\kappa_\lambda,\eta_\lambda,\chi_\lambda)$, for all $0\le\mu<\lambda<\infty$.
    \item[\rm(ii)] The generating functions $\varphi_q$ of $q$ and $\varphi_\eta$ of $\eta$ satisfy $\varphi_q(s)=s+\widetilde{\psi}(1-s)$ for some 
      $\widetilde{\psi}$ of the form {\rm(\ref{lkbm})} and $\varphi_\eta(s)=1-\widetilde{\phi}(1-s)$ for some $\widetilde{\phi}$ of the form 
      {\rm(\ref{lkim})}; $\kappa$ is geometrically divisible for all $\alpha>1/\widetilde{\psi}^\prime(\infty)$ if $\widetilde{\psi}^\prime(\infty)<\infty$,
      or for all $\alpha>0$ if $\widetilde{\psi}^\prime(\infty)=\infty$.
      Moreover, $\chi$ is also geometrically divisible
	  \begin{itemize}\item for all $\alpha>1/\widetilde{\phi}(\infty)$ if $\widetilde{\phi}(\infty)<\infty$;
        \item for all $\alpha>0$ if $\widetilde{\phi}(\infty)=\infty$.
      \end{itemize}
  \end{enumerate}
  In the setting of {\rm (i)} and {\rm (ii)}, a consistent family $(F_\lambda)_{\lambda\ge 0}$ of
  ${\rm GWI}(q_\lambda,\kappa_\lambda,\eta_\lambda,\chi_\lambda)$-forests can
  be constructed such that $(F_\mu,F_\lambda)\overset{{\rm (d)}}{=} (F_\lambda^{(1-\mu/\lambda)-\rm rdc},F_\lambda)$ for all
  $0\le\mu<\lambda<\infty$. 
\end{coro}

\noindent We omit the proof which is a straightforward combination of the proofs of Theorems \ref{thm2} and \ref{thm3}.

Similarly, the convergence results of Proposition \ref{prop4a} and \ref{prop4} find their analogue in this setting:

\begin{coro}\label{corolast2} Let $(Y_t^\lambda,t\ge 0)$ be the population size process in the setting of Corollary \ref{corolast}. Then
$$ \frac{Y_t^\lambda}{\psi^{-1}(\lambda)}\rightarrow Y_t\qquad\mbox{in distribution as $\lambda\rightarrow\infty$, for all $t\ge 0$,}
$$
where $(Y_t,t\ge 0)$ is a ${\rm CBI}(K,\psi,\widehat{K},\phi)$ with $Y_0\!=\!0$, for 
\begin{itemize}\item branching mechanism $\psi$ a linear transformations of $\widetilde{\psi}$ as in (\ref{psipsitilde})
  \item $K=(K_s,s\ge 0)$ with $\inf\{s\ge 0\colon K_s>V_\lambda\}\sim\kappa_\lambda$ for 
$V_\lambda\sim{\rm Exp}(c_\lambda)$ with $c_\lambda$ of Theorem \ref{thm3},
  \item immigration mechanism $\phi$ a linear transformation of $\widetilde{\phi}$ as in (\ref{full immi gen by phi psi})
  \item $\widehat{K}=(\widehat{K}_s,s\ge 0)$ with $\inf\{s\ge 0\colon \widehat{K}_s>\widehat{V}_\lambda\}\sim\chi_\lambda$ for $\widehat{V}_\lambda\sim{\rm Exp}(\paramgamma_\lambda)$ with $\paramgamma_\lambda$ of Theorem \ref{thm3}.
\end{itemize}
  If furthermore $\psi^\prime(0)>-\infty$ and $\phi^\prime(0)<\infty$, then the convergence holds in the almost sure sense.
\end{coro}

\noindent Like a ${\rm CBI}(\psi,\phi)$, a ${\rm CBI}(K,\psi,\widehat{K},\phi)$ can be constructed from a Poisson point process $(E^s,s\ge 0)$ in $\bD([0,\infty),[0,\infty))$
with intensity measure $d\Theta_{K,\psi}+\int_{(0,\infty)}\bP^x_{K,\psi}\Lambda(dx)$, where $(d,\Lambda)$ are the characteristics of $\phi$ in (\ref{lkim}). Also
consider an independent subordinator $\widehat{\sigma}$ or increasing random walk, in the $\phi(\infty)<\infty$ case, in fact $\widehat{K}_s=\inf\{t\ge 0\colon \widehat{\sigma}(t)>s\}$ clarifies the notation and the meaning of $\widehat{K}$. Then
$$Y_t=\sum_{s\le\widehat{K}_t}E_{t-\widehat{\sigma}(s)}^s$$
is a ${\rm CBI}(K,\psi,\hat{K},\phi)$ with $Y_0=0$. Like a ${\rm CSBP}(K,\psi)$, a ${\rm CBI}(K,\psi,\widehat{K},\phi)$ is non-Markovian, but admits a Markovian representation $(M,\vartheta)$ 
that records residual lifetimes as well as residual times to the next immigration, with values in $\bM([0,\infty))\times[0,\infty)$, where
$\vartheta_t=\widehat{\sigma}(\widehat{K}_t)-t$, so that
\beq&&\hspace{-0.5cm}\bE\left(\left.\!\exp\left\{\!-\!\!\int_{[0,\infty)}\!\!f(z)M_t(dz)-r\vartheta_t\!\right\}\right|M_0\!=\!m\!,\vartheta_0=s\right)\\
	&&=\bE\Bigg(\exp\Bigg\{-\int_{[0,t]}u_{t-z}(f)m(dz)-\int_{(t,\infty)}f(z-t)m(dz)\\
	&&\hspace{2.1cm}-\int_{[0,t-s]}\phi(u_{t-s-z}(f))d\widehat{K}_z-r\left(\widehat{\sigma}(\widehat{K}_{(t-s)^+})-(t-s)\right)\Bigg\}\Bigg),\eeq
where $u_t(f)$ is the unique nonnegative solution of (\ref{measint}). Then the process $Y_t=M_t([0,\infty))$ is a ${\rm CBI}(K,\psi,\hat{K},\phi)$. We leave any
further details including the proof of Corollary \ref{corolast2} to the reader.

%------------------------------------------------------Section 5-----------------------------------------------------------------------------

\bibliographystyle{abbrv}
\bibliography{CaW}
\end{document}